\documentclass[a4paper,12pt,oneside]{article}
\usepackage{main}

\begin{document}

\title{%
Gradient-based collective and structured nonlinear dimension reduction with quadratic surrogates%
\thanks{This project is funded by the ANR-DFG project COFNET (ANR-21-CE46-0015). This work was partially conducted within the France 2030 framework programme, Centre Henri Lebesgue ANR-11-LABX-0020-01.}}
\author{%
  A. Pasco${}^{1,\dagger}$,
  A. Nouy${}^{1}$}
\date{\medskip%
  \small %
  ${}^1$  \'Ecole Centrale de Nantes, Nantes Universit\'e, \\ Laboratoire de Mathématiques Jean Leray UMR CNRS 6629\\
  ${}^\dagger$ \emph{Present address:} Chair of High Performance Numerical Algorithms and Simulations, Institute of Mathematics, \'Ecole polytechnique f\'ed\'erale de Lausanne (EPFL), Lausanne, Switzerland\\
  \texttt{alexandre.pasco1702@gmail.com};
  \texttt{anthony.nouy@ec-nantes.fr}
}


\maketitle

\begin{abstract}

This paper is concerned with the approximation of continuously differentiable functions with high-dimensional input by a composition of two functions: a feature map that extracts few features from the input space, and a profile function that approximates the target function taking the features as its low-dimensional input.
We focus on the construction of structured nonlinear feature maps, that extract features on separate groups of variables, using a recently introduced gradient-based method that leverages Poincar\'e inequalities on nonlinear manifolds.
This method consists in minimizing a non-convex loss functional, which can be a challenging task, especially for small training samples.
We first investigate a collective setting, in which we construct a feature map suitable to a parametrized family of high-dimensional functions.
In this setting we introduce a new quadratic surrogate to the non-convex loss function and show an upper bound on the latter.
We then investigate a grouped setting, in which we construct separate feature maps for separate groups of inputs, and we show that this setting is almost equivalent to multiple collective settings, one for each group of variables.

\end{abstract}



\clearpage


\section{Introduction}
\label{sec:structured features introduction}

Recent decades have seen the development of increasingly accurate numerical models, but these are also increasingly costly to simulate.
However, for many purposes such as inverse problems, uncertainty quantification, or optimal design, many evaluations of these models are required.
A common approach is to use surrogate models instead, which aim to approximate the original model well while being cheap to evaluate.
Classical approximation methods, such as polynomials, splines, or wavelets, often perform poorly when the input dimension of the model is large, especially when few samples of the model are available.
Dimension reduction methods can help solve this problem.

This paper is concerned with two different settings in high-dimensional approximation.
Firstly, we consider a \emph{collective} dimension reduction setting, in which we aim to approximate functions from a parametrized family of continuously differentiable functions $u(\cdot, y) : \calX \rightarrow \bbR$ parametrized by some $y\in\calY$, where $\calX \subset \bbR^d$, $d\gg1$.  
We consider an approximation of the form
\[
    \hat u (\bfX, Y) = f(g(\bfX), Y),
\]
for some \emph{feature map} $g : \calX \rightarrow \bbR^m$, $m\ll d$, and a \emph{profile function} $f :\bbR^m \times \calY \rightarrow \bbR$, where the error is assessed in the $L^2(\calX \times \calY, \mu_{\calX} \otimes \mu_{\calY})$-norm for some probability distributions $\mu_{\calX}$ of $\bfX$ on $\calX$ and $\mu_{\calY}$ of $Y$ on $\calY$.
Secondly, we consider a \emph{grouped} or \emph{separated} dimension reduction setting, in which we aim to approximate a continuously differentiable function $u : \calX \rightarrow \bbR$ by splitting the input variables into $N$ groups, for some partition $S = \{\alpha_1, \cdots, \alpha_N\}$ of $\{1, \cdots, d\}$ containing disjoint multi-indices $\alpha_i \subset \{1, \cdots, d\}$, writing $\bfx = (\bfx_{\alpha})_{\alpha\in S}$ and $\calX = \times_{\alpha \in S} \calX_{\alpha}$.
We then consider an approximation of the form
\[
    \hat u (\bfX) 
    = f(g^{\alpha_1}(\bfX_{\alpha_1}), \cdots, g^{\alpha_N}(\bfX_{\alpha_N})),
\]  
for some feature maps $g^{\alpha} : \calX_{\alpha} \rightarrow \bbR^{m_{\alpha}}$ and some profile function $f : \times_{\alpha \in S} \bbR^{m_{\alpha}} \rightarrow \bbR$, where the error is assessed in the $L^2(\calX, \otimes_{\alpha \in S} \mu_{\alpha})$-norm for some probability distributions $\mu_{\alpha}$ of $\bfX_{\alpha}$ on $\calX_{\alpha}$, for all $\alpha \in S$.

Both the collective and the grouped settings can be seen as special cases of a more general dimension reduction setting $\hat u = f \circ g$, where a specific structure is imposed on the feature map.
Such structure may arise naturally from the original model, and allows for the incorporation of a priori knowledge in the feature map.

When the feature map is linear, i.e. $g(\bfx) = G^T \bfx$ for some $G \in \bbR^{d\times m}$, then $\hat u$ is a so-called \emph{ridge function} \cite{pinkusRidgeFunctions2015}, for which a wide range of methods have been developed.
The most classical one is the principal component analysis \cite{pearsonLinesPlanesClosest1901,hotellingAnalysisComplexStatistical1933}, with its grouped variant \cite{takaneComponentAnalysisDifferent1995}, which consists of choosing a $G$ that spans the dominant eigenspace of the covariance matrix of $\bfX$, without using information on $u$ itself.
Other statistical methods consist of choosing a $G$ that spans the \emph{central subspace}, such that $u(\bfX)$ and $\bfX$ are independent conditionally on $G^T \bfX$, which can be written in terms of the conditional measures $\mu_{(u(\bfX), \bfX) | G^T \bfX} = \mu_{u(\bfX) | G^T \bfX} \otimes \mu_{\bfX | G^T \bfX}$ almost surely.
Such methods are called \emph{sufficient dimension reduction} methods.
Major examples include \cite{liSlicedInverseRegression1991,denniscookMethodDimensionReduction2000,liDirectionalRegressionDimension2007}, and grouped variants include \cite{liGroupwiseDimensionReduction2010,guoGroupwiseDimensionReduction2015,liuStructuredOrdinaryLeast2017}.
We refer to \cite{liSufficientDimensionReduction2018} for a broad overview on sufficient dimension reduction.
Note that the collective setting can be seen as a special case of \cite{virtaSlicedInverseRegression2022}.

One problem with such methods is that they do not provide certification on the error one makes by approximating $u$ by a function of $G^T \bfx$.
Such certification can be obtained by leveraging Poincar\'e inequalities and gradient evaluations, leading to a bound of the form
\begin{equation}
\label{equ:upper bound approx error linear}
    \min_{\substack{f : \bbR^m \rightarrow \bbR \\ f \text{ measurable}}}
    \Expe{|u(\bfX) - f(G^T \bfX)|^2}
    \leq 
    C \Expe{\|\nabla u(\bfX)\|_2^2 - \|\Pi_{G} \nabla u(\bfX)\|_2^2},
\end{equation}
where $C>0$ depends on the distribution of $\bfX$, and where $\Pi_{G} := \Pi_{\spanv{G}} \in \bbR^{d\times d}$ denotes the orthogonal projector onto the column span of $G$.
The so-called \emph{active-subspace method} \cite{constantineActiveSubspaceMethods2014,zahmGradientBasedDimensionReduction2020} then consists of choosing a $G\in\bbR^{d\times m}$ that minimizes the right-hand side of the above equation, which turns out to be any matrix whose columns span the dominant eigenspace of $\Expe{\nabla u(\bfX) \nabla u(\bfX)^T}$.

Despite the theoretical and practical advantages of linear dimension reduction, some functions cannot be efficiently approximated with few linear features, for example $u(\bfx) = h(\|\bfx\|_2^2)$ for some $h\in\calC^1$.
For this reason, it may be worthwhile to consider nonlinear feature maps $g$.
Most aforementioned methods have been extended to nonlinear features, starting with the kernel principal component analysis \cite{scholkopfNonlinearComponentAnalysis1998}.
Nonlinear sufficient dimension reduction methods have also been proposed \cite{yi-renyehNonlinearDimensionReduction2009,leeGeneralTheoryNonlinear2013,wangFunctionalContourRegression2013,liNonlinearSufficientDimension2017,liDimensionReductionFunctional2022}, where the collective setting can again be seen as a special case of \cite{dongFrechetKernelSliced2022,yingFrechetSufficientDimension2022,zhangNonlinearSufficientDimension2024,zhangDimensionReductionFrechet2024}.
Gradient-based nonlinear dimension reduction methods have also been introduced, leveraging Poincar\'e inequalities \cite{bigoniNonlinearDimensionReduction2022,romorLocalApproachParameter2024,verdiereDiffeomorphismbasedFeatureLearning2025,nouySurrogatePoincareInequalities2026}, or not \cite{bridgesActiveManifoldsNonlinear2019,zhangLearningNonlinearLevel2019,gruberNonlinearLevelSet2021,romorKernelbasedActiveSubspaces2022,tengLevelSetLearning2023}.
In particular, an extension of \eqref{equ:upper bound approx error linear} to nonlinear feature maps was proposed in \cite{bigoniNonlinearDimensionReduction2022},
\begin{equation}
\label{equ:upper bound approx error nonlinear}
    \min_{\substack{f : \bbR^m \rightarrow \bbR \\ f \text{ measurable}}}
    \Expe{|u(\bfX) - f(g(\bfX))|^2}
    \leq 
    C \Expe{\|\nabla u(\bfX)\|_2^2 - \|\Pi_{\nabla g(\bfX)} \nabla u(\bfX)\|_2^2}
    := C \calJ(g),
\end{equation}
where $C>0$ depends on the distribution of $\bfX$ and the set of available feature maps, and where $\nabla g(\bfX) := (\nabla g_1(\bfX), \cdots \nabla g_m(\bfX)) \in \bbR^{d\times m}$ is the transposed jacobian matrix of $g$.
One issue in the nonlinear setting is that minimizing $\calJ$ over a set of nonlinear feature maps can be challenging as it is non-convex.
Circumventing this issue was the main motivation for \cite{nouySurrogatePoincareInequalities2026}, where quadratic surrogates to $\calJ$ were introduced and analyzed for a class of feature maps including polynomials.
The main contribution of the present work is to extend this approach to the collective setting.

In addition to the previously cited works, we refer the reader to \Cref{apx:poincare inequalities for dimension reduction} for discussions on the theoretical and practical aspects underlying \eqref{equ:upper bound approx error linear} and \eqref{equ:upper bound approx error nonlinear}.

Let us emphasize that the approaches described in this section are two-step procedures.
The feature map $g$ is learnt in a first step, without taking into account the class of profile functions used in the second step.
The second step consists of using classical regression tools to approximate $u$ as a function of $g(\bfx)$.
Alternatively, one may consider learning $f$ and $g$ simultaneously as in \cite{hokansonDataDrivenPolynomialRidge2018,lataniotisExtendingClassicalSurrogate2020}.

\subsection{Contributions and outline}

The first main contribution of the present work concerns the collective dimension reduction setting from \Cref{sec:collective dimension reduction}.
Applying the approach from \cite{bigoniNonlinearDimensionReduction2022} to $u(\cdot, y)$ for all $y\in\calY$ yields a collective variant of \eqref{equ:upper bound approx error nonlinear} with
\begin{equation}
\label{equ:def of J collective}
    \calJ_{\calX}(g) := \Expe{\|\nabla_{\bfx} u(\bfX,Y)\|_2^2 - \|\Pi_{\nabla g(\bfX)} \nabla_{\bfx} u(\bfX,Y) \|_2^2},
\end{equation}
which is again a non-convex function for nonlinear feature maps.
Following \cite{nouySurrogatePoincareInequalities2026}, we introduce a new quadratic surrogate in order to circumvent this problem,
\begin{equation}
\label{equ:def L collective}
    \calL_{\calX,m}(g) := \Expe{
        \lambda_1(M(\bfX)) (\|\nabla g(\bfX)\|_F^2 - \| \Pi_{V_m(\bfX)} \nabla g(\bfX)\|_F^2 )
    },
\end{equation}
where the columns of $V_m(\bfX)\in\bbR^{d\times m}$ are the $m$ principal eigenvectors of the conditional uncentered second moment matrix $M(\bfX) = \Expe[Y]{\nabla_{\bfx} u(\bfX, Y) \nabla_{\bfx} u(\bfX, Y)^T} \in \bbR^{d\times d}$, with $\lambda_1(M(\bfX))$ its largest eigenvalue.
We show that for non-constant polynomial feature maps of degree at most $\ell+1$,
\[
    0 \leq \calJ_{\calX}(g) - \varepsilon_m 
    \lesssim \calL_{\calX, m}(g)^{\frac{1}{1+ 2\ell m}},
\]
where $\varepsilon_m = \sum_{i=m+1}^d\Expe{\lambda_i(M(\bfX))}$ is a lower bound on $\calJ_{\calX}$ that does not depend on the feature maps.
We then show that if $g(\bfx) = G^T \Phi(\bfx)$ for some $\Phi \in \calC^1(\calX, \bbR^K)$ and some $G\in\bbR^{K\times m}$ then
\[
\begin{aligned}
    \calL_{\calX, m}(g) 
    & = \Trace{G^T H_{\calX, m} G}, 
    \\
    H_{\calX, m} 
    &= \Expe{
        \lambda_1(M(\bfX)) 
        \nabla \Phi(\bfX)^T 
        \big(
            I_d - V_m(\bfX) V_m(\bfX)^T
        \big)
        \nabla \Phi(\bfX)
    } \in \bbR^{K\times K},
\end{aligned}
\]
which means that minimizing $\calL_{\calX, m}$ is equivalent to finding the eigenvectors associated to the smallest eigenvalues of $H_{\calX ,m}$.
There are three main differences with the surrogate-based approach from \cite{nouySurrogatePoincareInequalities2026}.
Firstly, estimating $V_m(\bfX)$ and $\lambda_1(M(\bfX))$ requires a sample of $Y$ for each sample point of $X$, using for example a tensorized sample of the form $(\bfx^{(i)}, y^{(j)})_{1\leq i\leq n_{\calX}, 1\leq j\leq n_{\calY}}$ with size $n=n_{\calX}n_{\calY}$, which may be prohibitive and is the main limitation of our approach.
Secondly, the collective setting allows for richer information on $u$ for fixed $\bfx$, so that the surrogate $\calL_{\calX,m}$ can be directly used in the case $m>1$, while \cite{nouySurrogatePoincareInequalities2026} relies on successive surrogates to learn one feature at a time.
Thirdly, we only show that our new surrogate can be used as an upper bound, while \cite{nouySurrogatePoincareInequalities2026} provided both lower and upper bounds.

The second main contribution concerns near-optimality results for the grouped dimension reduction setting, presented in \Cref{sec:two variables approach,sec:multiple variables approach}.
By drawing a parallel with tensor approximation, more precisely with the higher order singular value decomposition (HOSVD), we show that both grouped dimension-reduction problems can be decomposed, up to near-equivalence, into multiple collective settings.

The rest of this paper is organized as follows.
First in \Cref{sec:collective dimension reduction} we introduce and analyze our new quadratic surrogate for collective dimension reduction.
Then in \Cref{sec:two variables approach} and \Cref{sec:multiple variables approach}, we investigate grouped settings with two groups and more groups of variables, respectively, and show that they are nearly equivalent to multiple collective dimension reduction settings.
Then in \Cref{sec:toward hierarchical formats} we briefly discuss extensions toward hierarchical formats, although we only provide pessimistic examples.
Then in \Cref{sec:numerical experiments structured dim red} we illustrate the collective dimension reduction setting on a numerical example.
Finally, in \Cref{sec:conclusion structured features} we summarize the analysis and observations and we discuss perspectives.


\section{Collective dimension reduction}
\label{sec:collective dimension reduction}

In this section, we consider a dimension reduction problem for $u : \calX \times \calY \rightarrow \bbR$ with respect to the first variable $\bfX$, in order to approximate $u(\bfX, Y)$ in the space  $L^2(\calX \times \calY, \mu_{\calX} \otimes \mu_{\calY})$.
We want this dimension reduction to be collective, in the sense that the feature maps for $\bfX$ shall be the same for any realization of the random function $u_Y := u(\cdot, Y)$.
In other words, we consider an approximation $\hat u_y : \calX \rightarrow \bbR$ of the form
\[
    \hat u_y : \bfx \mapsto 
    f(g(\bfx), y)
\]
with $g : \calX \rightarrow \bbR^m$ and $f : \bbR^m \times \calY \rightarrow \bbR$ belonging respectively to some classes of feature maps $\calG_m \subset \calC^1(\calX, \bbR^m)$ and profile functions $\calF_m$.
Following the approach from \cite{bigoniNonlinearDimensionReduction2022}, we consider no restriction besides measurability on the profile functions, so that we would want to construct a feature map that minimizes
\begin{equation}
\label{equ:def reconstruction error collective}
    \calE_{\calX}(g) := 
    \min_{\substack{f: \bbR^m \times \calY \rightarrow \bbR \\ f \text{ measurable}}}
    \Expe{|u_Y(\bfX) - f(g(\bfX), Y)|^2},
\end{equation}
where the minimum in the above equation is obtained by the conditional expectation $f_g : (\bfz, y) \mapsto \Expe{u_Y(\bfX) | (g(\bfX), Y) = (\bfz,y)}$.
Now, under suitable assumptions on $\calG_m$, we can apply \cite[Proposition 2.9]{bigoniNonlinearDimensionReduction2022} on $u_Y$ and take the expectation over $Y$ to obtain
\begin{equation} 
\label{equ:E bounded by J collective}
    \calE_{\calX}(g)
    \leq C \calJ_{\calX}(g)
    = C \Expe{\|\nabla u_Y(\bfX)\|_2^2 - \|\Pi_{\nabla g(\bfX)} \nabla u_Y(\bfX) \|_2^2},
\end{equation}
with the same $C$ as in \eqref{equ:upper bound approx error nonlinear}.
We refer to Appendix \ref{apx:poincare inequalities for dimension reduction} for a discussion on this constant.
Note that we can also write $\calJ_{\calX}$ as $\calJ_{\calX}(g) = \calJ(\tilde g)$ with $\calJ$ defined in \eqref{equ:upper bound approx error nonlinear} and $\tilde g : (\bfx, y) \mapsto (g(\bfx), y)$.

In the rest of this section, we design a quadratic surrogate to $\calJ_{\calX}$ in a manner similar to \cite{nouySurrogatePoincareInequalities2026}.
Firstly, in \Cref{subsec:truncation in the poincare inequality-based loss} we introduce a truncated version $\calJ_{\calX,m}$ of $\calJ_{\calX}$, and we show that minimizing $\calJ_{\calX,m}$ is almost equivalent to minimizing $\calJ_{\calX}$.
Secondly, in \Cref{subsec:quadratic surrogate to the truncated loss} we introduce a new quadratic function $\calL_{\calX, m}$ as a surrogate to $\calJ_{\calX, m}$, and we show that it can be used to upper bound $\calJ_{\calX, m}$ for bi-Lipschitz or polynomial feature maps.
Thirdly, in \Cref{subsec:minimizing the surrogate collective} we show that, when the feature map's coordinates are taken as orthonormal elements of some finite dimensional vector space of functions, then minimizing $\calL_{\calX, m}$ is equivalent to solving a generalized eigenvalue problem.

\begin{rmrk}
    A particular case of the collective setting is the vector valued setting.
    Indeed, approximating $v: \bfx \mapsto (v_1(\bfx), \cdots, v_n(\bfx)) \in \bbR^n$ in $L^2(\calX, \mu_{\calX}; \bbR^n)$ is equivalent to approximating $u: (\bfx, y) \mapsto v_y(\bfx)$ in $L^2(\calX \times \calY, \mu_{\calX}\otimes \mu_{\calY})$ with $\mu_{\calY}$ the uniform measure on $\calY = \{1, \cdots, n \}$.
\end{rmrk}

\begin{rmrk}
    In this section we assume that $\mu_{\calY}$ is a probability measure, which allowed us to stay in a rather classical setting and to simplify notations.
    However, this assumption is most probably not necessary, as one should be able to derive the same analysis with a more general measure $\mu_{\calY}$, although it would require some rewriting.
    We leave this aspect to future investigation.
\end{rmrk}

\subsection{Truncation of the Poincar\'e inequality based loss}
\label{subsec:truncation in the poincare inequality-based loss}

In this section, we introduce a truncated version $\calJ_{\calX, m}$ of $\calJ_{\calX}$ defined in \eqref{equ:def of J collective}, and we show that minimizing this truncated version is almost equivalent to minimizing $\calJ_{\calX}$.
The intuition behind this truncation step is that the collective setting has an intrinsic limit in the dimension reducibility, compared to the non-collective setting (that is, when $Y$ is constant). 
Indeed, in the latter there always exists an optimal scalar feature map $g\in\calC^{1}(\calX, \bbR^m)$ such that $\calJ(g) = 0$, for example consider $g=u$.
This may not be the case in the collective setting, as demonstrated in Example \ref{exa:dimension reduction limit}.
\begin{xmpl}[Two features required]
\label{exa:dimension reduction limit}
Consider a function of the form $u_y(\bfx) = u^0(\bfx) + y u^1(\bfx)$ with two non constant functions $u^0$ and $u^1$.
Obviously, taking $m=2$ with $g(\bfx) = (u^0(\bfx), u^1(\bfx))$ and $f(\bfz,y)=z_1 + yz_2$ is sufficient to obtain $u_y(\bfx)=f(g(\bfx),y)$.
This cannot be true for $m=1$.
Indeed, we have $\nabla u_y(\bfx) = \nabla u^0(\bfx) + y \nabla u^1(\bfx)$, thus
\[
    \nabla u_y(\bfX) 
    \in \spanv{\nabla u_y(\bfX), y\in \calY}
    = \spanv{\nabla u^0(\bfX), \nabla u^1(\bfX)}.
\]
Hence, if $\nabla u^0(\bfX)$ and $\nabla u^1(\bfX)$ are not almost surely colinear, then $\spanv{\nabla u_y(\bfX), y\in \calY}$ has dimension two.
However, if we could write $u_y(\bfx)=f(g(\bfx),y)$ for some $g$ and $f$ with $m=1$, then since $\nabla u_y(\bfx) \in \spanv{\nabla g(\bfx)}$, the dimension of $\spanv{\nabla u_y(\bfX), y\in \calY}$ would be one, which contradicts the previous statement.
Hence, the function $u$ in this example requires (at least) two features in $\bfX$ to be exactly represented.
\end{xmpl}

The first step is then to investigate a lower bound on $\calJ_{\calX}$ that does not depend on the feature maps considered.
This can be obtained by searching for a matrix $V_m(\bfX)$ whose column span is better than any column span of $\nabla g(\bfX)$ for any possible $g\in\calC^{1}(\calX, \bbR^m)$.
In Example \ref{exa:dimension reduction limit} one could choose any $V_2(\bfX)$ whose columns span $\spanv{\nabla u^0(\bfX), \nabla u^1(\bfX)}$.
We thus naively define $V_m(\bfX)\in\bbR^{d\times m}$ as a matrix satisfying 
\begin{equation}
\label{equ:def princial components gradient u}
    \Expe[Y]{\|\Piperp_{V_m(\bfx)} \nabla u_Y(\bfx)\|_2^2}
    = 
    \min_{W \in \bbR^{d\times m}}
    \Expe[Y]{\|\Piperp_{W} \nabla u_Y(\bfx)\|_2^2},
\end{equation}
where $\Piperp_{\nabla g(\bfx)} := I_d - \Pi_{\nabla g(\bfx)}$.
By definition, $\Expe{\|\Piperp_{V_m(\bfX)} \nabla u_Y(\bfX)\|_2^2} \leq \calJ_{\calX}(g)$ for any $g\in\calC^{1}(\calX, \bbR^m)$.
It turns out that $V_m(\bfx)$ is commonly known as the principal components matrix of $\nabla u_Y(\bfx)$, and can be defined as $V_m(\bfx) = (v^{(1)}(\bfx) , \cdots , v^{(m)}(\bfx))$ where $(v^{(i)}(\bfx))_{1\leq i\leq d} \subset \bbR^d$ are the eigenvectors associated to $\lambda_1(M(\bfx)) \geq \cdots \geq  \lambda_m(M(\bfx)) \geq 0$, the eigenvalues of the symmetric positive semidefinite uncentered second moment matrix 
\begin{equation}
\label{equ:def covariance matrix gradient u}
    M(\bfx) :=  
    \Expe[Y]{
    \nabla u_Y(\bfx) 
    \nabla u_Y(\bfx)^T
    } 
    \in \bbR^{d\times d}.
\end{equation}
Since $M(\bfx)$ is real symmetric, we assume that its eigenvectors $(v^{(1)}(\bfx) , \cdots , v^{(d)}(\bfx))$ are orthonormal, which implies $V_m(\bfx)^T V_m(\bfx) = I_m$.
By the properties of the eigenvectors, taking the expectation over $\bfX$ yields the following lower bound on $\calJ_{\calX}(g)$ in terms of the eigenvalues of the above matrix,
\begin{equation}
\label{equ:J collective lower bound eigen values}
\begin{aligned}
    \varepsilon_m
    &:= \Expe{\|\Piperp_{V_m(\bfX)} \nabla u_Y(\bfX)\|_2^2}
    = \Trace{\Expe{\Piperp_{V_m(\bfX)} \nabla u_Y(\bfX) \nabla u_Y(\bfX)^T \Piperp_{V_m(\bfX)}}} \\
    & = \Trace{\Expe{\Piperp_{V_m(\bfX)} M(\bfX) \Piperp_{V_m(\bfX)}}}
    = \sum_{i=m+1}^{d} \Expe{\lambda_i(M(\bfX))}
    \leq \calJ_{\calX}(g).
\end{aligned}
\end{equation}
Note that in the first line of the above equation we used the property that $\Expe{\|\bfZ\|^2_2} = \Trace{\Expe{\bfZ \bfZ^T}}$ for any square-integrable random vector $\bfZ$.
We denote this property as the \textit{trace identity for the second moment} as it will be used several times in the present work.
Note also that we further discuss the computation of $M(\bfX)$ and $V_m(\bfX)$, which is the major computational aspect of our approach, at the end of \Cref{subsec:minimizing the surrogate collective}.
Also, we refer to Figure \ref{fig:illustration dim red collective} for an illustrated example.

\begin{rmrk}
The quantity $\varepsilon_m$ serves as an indicator of the dimension reducibility of the collective dimension reduction problem.
In Example \ref{exa:dimension reduction limit} we have $\varepsilon_1 > 0$, which implies that the problem cannot be reduced to a single feature ($m=1$).
On the other hand, the problem can be reduced to $m=2$ features, which implies $\varepsilon_2 = 0$.
Note that for the converse statement to hold, one would need to show that there exists $g\in\calC^{1}(\calX, \bbR^m)$ such that $\spanv{\nabla g(\bfX)} = \spanv{V_m(\bfX)}$ almost surely while satisfying a connectivity assumption on its level set.
See \cite[Proposition 2.3]{bigoniNonlinearDimensionReduction2022} for the connectivity assumption.
\end{rmrk}

We thus propose to build some feature map $g\in\calG_m$ whose gradient is aligned with $\Pi_{V_m(\bfX)}\nabla u_Y(\bfX)$ instead of $\nabla u_Y(\bfX)$, by defining the truncated version of $\calJ_{\calX}$ as
\begin{equation}
\label{equ:def J collective truncated}
    \calJ_{\calX, m}(g) :=
    \Expe{
        \|\Piperp_{\nabla g(\bfX)} \Pi_{V_m(\bfX)} \nabla u_Y(\bfX)\|_2^2
    }.
\end{equation}
The first interesting property of $\calJ_{\calX,m}$ is that it is almost equivalent to $\calJ_{\calX}$ as a measure of quality of a feature map $g\in\calG_m$.
In particular, any minimizer of $\calJ_{\calX, m}$ is almost a minimizer of $\calJ_{\calX}$.
These properties are stated in Proposition \ref{prop:truncated J almost optimal}.

\begin{prpstn}
\label{prop:truncated J almost optimal}
Let $\calJ_{\calX}$, $\calJ_{\calX, m}$ and $\varepsilon_m$ be as defined respectively in \eqref{equ:def of J collective}, \eqref{equ:def J collective truncated} and \eqref{equ:J collective lower bound eigen values}.
Then for any $g\in\calG_m$,
\begin{equation}
\label{equ:truncated J equivalent to J}
    \frac{1}{2} 
    (\calJ_{\calX, m}(g) + \varepsilon_m)
    \leq \calJ_{\calX}(g)
    \leq \calJ_{\calX, m}(g) 
    + \varepsilon_m.
\end{equation}
Moreover, if $g^*$ is a minimizer of $\calJ_{\calX, m}$ over $\calG_m$ then
\begin{equation}
\label{equ:truncated J almost optimal}
    \calJ_{\calX}(g^*) \leq 2 \inf_{g\in \calG_m} \calJ_{\calX}(g).
\end{equation}
\end{prpstn}
\begin{proof}
    Firstly, for all $\bfx \in \calX$, since the eigenvectors of $M(\bfx)$ are orthogonal, we have that
    \[
        \Expe[Y]{(\Pi_{V_m(\bfx)} \nabla u_Y(\bfx))(\Piperp_{V_m(\bfx)} \nabla u_Y(\bfx))^T} =
        \Pi_{V_m(\bfx)} \Expe[Y]{\nabla u_Y(\bfx) \nabla u_Y(\bfx)^T} \Piperp_{V_m(\bfx)}
        = 0.
    \]
    Then, by multiplying the first term in the above equation on the left and right by $\Piperp_{\nabla g(\bfX)}$, and by taking the trace of the result, we obtain that the functions 
    \[
    y \mapsto \Piperp_{\nabla g(\bfx)} \Pi_{V_m(\bfx)} \nabla u_y(\bfx) 
    \quad \text{and} \quad 
    y \mapsto \Piperp_{\nabla g(\bfx)} \Piperp_{V_m(\bfx)} \nabla u_y(\bfx)
    \] 
    are orthogonal to each other in $L^2(\mathcal{Y}, \mu_{\mathcal{Y}}, \mathbb{R}^d)$.
    Since they sum to $y\mapsto \Piperp_{\nabla g(\bfx)} \nabla u_y(\bfx)$, the Pythagorean theorem yields
    \[
    \Expe[Y]{\|\Piperp_{\nabla g(\bfx)} \nabla u_Y(\bfx)\|_2^2}
    = \Expe[Y]{\|\Piperp_{\nabla g(\bfx)} \Pi_{V_m(\bfx)} \nabla u_Y(\bfx)\|_2^2}
    + \Expe[Y]{\|\Piperp_{\nabla g(\bfx)} \Piperp_{V_m(\bfx)} \nabla u_Y(\bfx)\|_2^2}.
    \]
    Then, by taking the expectation with respect to $\bfX$ in the above equation, we obtain
    \[
    \calJ_{\calX}(g) 
    = \calJ_{\calX, m}(g)
    + \Expe{
        \|\Piperp_{\nabla g(\bfX)} \Piperp_{V_m(\bfX)} 
        \nabla u_Y(\bfX)\|_2^2.
    }
    \]
    By observing that the second term on the right-hand side of the above equality is positive and upper bounded by $\varepsilon_m$ since $\|\Piperp_{\nabla g(\bfX)}\|_2 \leq 1$, we obtain
    \[
        \calJ_{\calX, m}(g) 
        \leq \calJ_{\calX}(g)
        \leq \calJ_{\calX, m}(g) + \varepsilon_m.
    \]
    Thus, summing the above inequalities with $\varepsilon_m \leq \calJ_{\calX}(g)$ from \eqref{equ:J collective lower bound eigen values} yields the desired inequality \eqref{equ:truncated J equivalent to J}.
    Finally, by using the right-hand inequality in \eqref{equ:truncated J equivalent to J}, the minimizing property of $g^*$, and the left-hand inequality in \eqref{equ:truncated J equivalent to J}, we obtain
    \[
        \calJ_{\calX}(g^*) 
        \leq \calJ_{\calX, m}(g^*) + \varepsilon_m
        \leq \calJ_{\calX, m}(g) + \varepsilon_m
        \leq 2\calJ_{\calX}(g),
    \]
    and taking the infimum over $g\in\calG_m$ yields the desired inequality \eqref{equ:truncated J almost optimal}.
\end{proof}

The second interesting property of $\calJ_{\calX, m}$ is that it is better suited to designing a quadratic surrogate using a similar approach to \cite{nouySurrogatePoincareInequalities2026}, which is the topic of the next \Cref{subsec:quadratic surrogate to the truncated loss}.

\subsection{Quadratic surrogate to the truncated loss}
\label{subsec:quadratic surrogate to the truncated loss}

In this section, inspired by \cite[Section 4]{nouySurrogatePoincareInequalities2026}, we detail the construction of a new quadratic surrogate which can be used to upper bound $\calJ_{\calX, m}$.
The first step toward this new surrogate is the following lemma.

\begin{lmm}
\label{lem:norm projection generalized}
    Let $n,m\leq d$ and let $V \in \bbR^{d\times n}$ and $W \in \bbR^{d\times m}$ be matrices such that $V^T V = I_n$ and $W^T W = I_m$.
    Then it holds
    \[
        \|\Piperp_W V\|_F^2 
        = \|\Piperp_V W\|_F^2 + (n-m).
    \]
\end{lmm}
\begin{proof}
    First, since $V$ is orthonormal we have that 
    $n = \|V\|_F^2 = \|\Piperp_W V\|_F^2 + \|\Pi_W V\|_F^2$.
    Similarly, it holds $m = \|W\|_F^2 = \|\Piperp_V W\|_F^2 + \|\Pi_V W\|_F^2$.
    Moreover, by assumption on $V$ and $W$ we have that $\Pi_V = V V^T$ and $\Pi_W = W W^T$, thus
    $\|\Pi_V W\|_F^2 = \|\Pi_W V\|_F^2 = \|V^T W\|_F^2$.
    Combining those two observations gives 
    $n - \|\Piperp_W V\|_F^2 = m - \|\Piperp_V W\|_F^2$, which yields the desired result.
\end{proof}

The case $n=m$ in the above Lemma \ref{lem:norm projection generalized} is of particular interest for us.
It is interpreted as an equality between the projection error of $V$ onto $\spanv{W}$ and the projection error of $W$ onto $\spanv{V}$.
We apply it to $W(\bfX) \in \bbR^{d\times m}$, whose column span is the same as the one of $\nabla g(\bfX)$, and $V_m(\bfX)\in\bbR^{d\times m}$ defined in \eqref{equ:def princial components gradient u}, which yields Lemma \ref{lem:bounds on J collective truncated} below.
The main problem is that the above interpretation is perturbed by the non orthogonality of $\nabla g(\bfX)$.

\begin{lmm}
\label{lem:bounds on J collective truncated}
    Let $g \in \calC^1(\calX, \bbR^m)$ such that $\mathrm{rank}(\nabla g(\bfX))=m$ almost surely.
    Then, with $M(\bfX)$ and $V_m(\bfX)$ as defined in \eqref{equ:def covariance matrix gradient u} and \eqref{equ:def princial components gradient u} respectively,
    \begin{equation}
    \begin{aligned}
        \calJ_{\calX, m}(g)
        \geq \Expe{
            \frac{\lambda_m(M(\bfX))}{\sigma_1(\nabla g(\bfX))^2}
            \|\Piperp_{V_m(\bfX)} \nabla g(\bfX) \|_F^2
        },
        \\
        \calJ_{\calX, m}(g)
        \leq \Expe{
            \frac{\lambda_1(M(\bfX))}{\sigma_m(\nabla g(\bfX))^2}
            \|\Piperp_{V_m(\bfX)} \nabla g(\bfX) \|_F^2
        }.
    \end{aligned}
    \end{equation}
\end{lmm}
\begin{proof}
    First, by using the trace identity for the second moment in the definition of $\calJ_{\calX, m}(g)$, then using the independence of $\bfX$ and $Y$, we obtain
    \[
    \begin{aligned}
    \calJ_{\calX, m}(g)
        &= \Expe{\Trace{
            \Piperp_{\nabla g(\bfX)}
            \Pi_{V_m(\bfX)}
            \nabla u_Y(\bfX)
            \nabla u_Y(\bfX)^T
            \Pi_{V_m(\bfX)}
            \Piperp_{\nabla g(\bfX)}
        }} \\
        &= \Expe[\bfX]{\Trace{
            \Piperp_{\nabla g(\bfX)}
            \Pi_{V_m(\bfX)}
            \Expe[Y]{
                \nabla u_Y(\bfX)
                \nabla u_Y(\bfX)^T
            }
            \Pi_{V_m(\bfX)}
            \Piperp_{\nabla g(\bfX)}
        }} \\
        &= \Expe[\bfX]{\Trace{
            \Piperp_{\nabla g(\bfX)}
            \Pi_{V_m(\bfX)}
            M(\bfX)
            \Pi_{V_m(\bfX)}
            \Piperp_{\nabla g(\bfX)}
        }}. \\
    \end{aligned}
    \]
    From now on in this proof, since there is no dependency on $Y$ anymore, let us not write the dependency in $\bfX$.
    Then, since $V_m$ is defined as the matrix whose columns are the dominant eigenvectors of $M$, we have $\Pi_{V_m} = V_m V_m^T$ and $V_m^T M V_m = D_m(M) := \mathrm{diag}((\lambda_i(M))_{1\leq i\leq m})$.
    Then, using the property of the trace of a product, then expressing the trace in the canonical basis $(e_i)_{1\leq i\leq m}$ of $\bbR^m$ with $(e_i)_j = \delta_{i,j}$, then using $D_m(M)e_i = \lambda_i(M)e_i$, we obtain
    \[
    \begin{aligned}
    \calJ_{\calX, m}(g)
        &= \Expe{\Trace{
            \Piperp_{\nabla g}
            V_m D_m(M) V_m^T
            \Piperp_{\nabla g}
        }}
        = \Expe{\Trace{
            D_m(M) V_m^T \Piperp_{\nabla g} V_m
        }} \\
        & = \sum_{1\leq i\leq m} \Expe{
            e_i^T D_m(M) V_m^T \Piperp_{\nabla g} V_m e_i
        }
        = \sum_{1\leq i\leq m} \Expe{
            \lambda_i(M) e_i^T V_m^T \Piperp_{\nabla g} V_m e_i
        } \\
        &= \sum_{1\leq i\leq m} \Expe{ 
            \lambda_i(M)
            v^{(i)~ T}\Piperp_{\nabla g} v^{(i)}
        }
        = \sum_{1\leq i\leq m} \Expe{ 
            \lambda_i(M)
            \|\Piperp_{\nabla g} v^{(i)}\|_2^2
        }.
    \end{aligned}
    \]
    Hence, denoting $v^{(i)}$ the $i$-th column of $V_m$, we obtain
    \begin{equation}
    \label{equ: expand Jx}
    \calJ_{\calX, m}(g) 
    = \sum_{1\leq i\leq m} \Expe{ 
        \lambda_i(M)
        \|\Piperp_{\nabla g} v^{(i)}\|_2^2
    }.
    \end{equation}
    Then, bounding the $m$ first eigenvalues of $M$ and identifying the squared Frobenius norm yields
    \begin{equation}
    \label{equ: bound Jx with V proj dg}
        \Expe{\lambda_m(M) \| \Piperp_{\nabla g} V_m \|_F^2}
        \leq \calJ_{\calX, m}(g) \leq 
        \Expe{\lambda_1(M) \| \Piperp_{\nabla g} V_m \|_F^2}.
    \end{equation}
    Let us now provide a lower and an upper bound of $\| \Piperp_{\nabla g} V_m \|_F^2$.
    Write the singular value decomposition of $\nabla g$ as $\nabla g = W \Lambda U^T$.
    Then, by applying Lemma \ref{lem:norm projection generalized} with $V=V_m$ and $W = \nabla g \Lambda^{-1} U$, which both have $m$ orthonormal columns, then by expressing the Frobenius norm in the canonical basis,
    we obtain
    \[
        \| \Piperp_{\nabla g} V_m \|_F^2
        = \| \Piperp_{ V_m} W \|_F^2
        = \| \Piperp_{ V_m} \nabla g U \Lambda^{-1}\|_F^2
        = \sum_{1\leq i \leq m} \sigma_i(\nabla g)^{-2} \| \Piperp_{ V_m} \nabla g U e_i\|_2^2.
    \]
    Then, by lower and upper bounding the singular values $\sigma_i(\nabla g)$, then by recognizing the expression of $\| \Piperp_{ V_m} \nabla g U\|_F^2$, which is equal to $\| \Piperp_{ V_m} \nabla g\|_F^2$ since $U U^T = I_m$ by definition, we obtain
    \begin{equation}
    \label{equ: bound V proj dg by dg proj V}
        \sigma_1(\nabla g)^{-2}
        \| \Piperp_{V_m} \nabla g\|_F^2
        \leq 
        \| \Piperp_{\nabla g} V_m \|_F^2
        \leq 
        \sigma_m(\nabla g)^{-2}
        \| \Piperp_{V_m} \nabla g\|_F^2,
    \end{equation}
    which combined with the previous inequalities on $\calJ_{\calX, m}(g)$ yields the desired result.
\end{proof}

In view of Lemma \ref{lem:bounds on J collective truncated}, we propose to define a new surrogate, with $M(\bfX)$ and $V_m(\bfX)$ defined in \eqref{equ:def covariance matrix gradient u} and \eqref{equ:def princial components gradient u} respectively, 
\begin{equation}
    \calL_{\calX,m}(g) :=
    \Expe{
        \lambda_1(M(\bfX))
        \|\Piperp_{V_m(\bfX)} \nabla g(\bfX)\|_F^2
    }.
\end{equation}
A first key property of this surrogate is that $g \mapsto \calL_{\calX,m}(g)$ is quadratic, and its minimization boils down to minimizing a generalized Rayleigh quotient when $g(\bfx)=G^T \Phi(\bfx)$ for some fixed $\Phi \in\calC^1(\calX, \bbR^K)$, $K\geq d$, as shown in \Cref{subsec:minimizing the surrogate collective}.
A second key property is that we can use $\calL_{\calX,m}$ to upper bound $\calJ_{\calX,m}$ for bi-Lipschitz or polynomial feature maps, as shown in \Cref{subsec:the surrogate as an upper bound}.
However, we are not able to provide the converse inequalities, meaning upper bounding $\calL_{\calX,m}$ with $\calJ_{\calX,m}$.

Finally, note that it remains consistent with the case $m=1$ from \cite[Section 4]{nouySurrogatePoincareInequalities2026}, as mentioned in Remark \ref{rmk:collective surrogate for one feature}.
Still, the current setting raises some additional questions, as pointed out in Remark \ref{rmk:not one choice for surrogate}.

\begin{rmrk}
\label{rmk:collective surrogate for one feature}
    Let us briefly show that Lemma \ref{lem:bounds on J collective truncated} and the new surrogate \eqref{equ:def L collective} remain consistent with the setting $m=1$ and $u_Y = u$ from \cite[Section 4]{nouySurrogatePoincareInequalities2026}.
    In this setting, we first observe that $\calJ(g) = \calJ_{\calX, 1}(g)$, and that the two inequalities in Lemma \ref{lem:bounds on J collective truncated} are actually equalities.
    Also, $\lambda_1(M(\bfX)) = \|\nabla u(\bfX)\|_2^2$ and $\sigma_1(\nabla g(\bfX)) = \|\nabla g(\bfX)\|_2$.
    As a result $\calL_{\calX,1}$ is exactly the surrogate from \cite[Section 4]{nouySurrogatePoincareInequalities2026},
    \[
        \calL_1(g)
        = \Expe{\|\nabla u(\bfX)\|_2^2 \|\Piperp_{\spanv{\nabla u(\bfX)}} \nabla g(\bfX) \|_2^2}.
    \]
\end{rmrk}

\begin{rmrk}
\label{rmk:not one choice for surrogate}
    A difference with the situation in \cite{nouySurrogatePoincareInequalities2026} is that there was somehow a natural choice of surrogate.
    This is not the case anymore, as one can legitimately replace $\lambda_1(M(\bfX))$ by any weighting $w(\bfX)$ such that $\lambda_m(M(\bfX)) \leq w(\bfX) \leq \lambda_1(M(\bfX))$.
    However, this choice will influence the available bounds, as choosing $\lambda_1(M(\bfX))$ allows to naturally obtain an upper bound on $\calJ_{\calX, m}$, while choosing $\lambda_m(M(\bfX))$ allows to naturally obtain a lower bound on $\calJ_{\calX, m}$.
    Since we want to minimize $\calJ_{\calX, m}(g)$, we have chosen the first option.
    Let us mention that one could obtain both upper and lower bounds if concentration inequalities on $\lambda_1(M(\bfX)) / \lambda_m(M(\bfX))$ were available.
\end{rmrk}

\subsection{The surrogate as an upper bound}
\label{subsec:the surrogate as an upper bound}

The intuition behind the surrogate $\calL_{\calX, m}$ is that instead of aligning $\spanv{\nabla g(\bfX)}$ with $\nabla u_Y(\bfX)$, which is the qualitative goal behind \eqref{equ:def of J collective} or \eqref{equ:def J collective truncated}, it aims to align the $m$ columns of $\nabla g(\bfX)$ with the principal eigenspace $\spanv{V_m(\bfX)}$.
In this section, we show that these two goals are linked to each other.
Let us first provide a result in the context of exact recovery, stated in Proposition \ref{prop:exact recovery collective} below.
This result applies to Example \ref{exa:dimension reduction limit}.

\begin{prpstn}
\label{prop:exact recovery collective}
    Assume that rank$(M(\bfX)) \geq m$ almost surely, with $M(\bfX)$ as defined in \eqref{equ:def covariance matrix gradient u}.
    Let $g\in\calC^1(\calX, \bbR^m)$ be such that $\mathrm{rank}(\nabla g(\bfX))=m$ almost surely.
    Then
    \[
        \calJ_{\calX, m}(g) = 0
        \iff
        \calL_{\calX, m}(g) = 0.
    \] 
\end{prpstn}
\begin{proof}
    Under the assumptions, we have that both $\lambda_m(M(\bfX))$ and $\sigma_m(\nabla g(\bfX))^2$ are almost surely strictly positive, so their ratio is almost surely finite and strictly positive.
    Then Lemma \ref{lem:bounds on J collective truncated} yields that $\calJ_{\calX, m}(g) = 0$ if and only if $\|\Piperp_{V_m(\bfX)} \nabla g(\bfX)\|_F^2 = 0$ almost surely.
    Finally, since $0 < \lambda_m(M(\bfX)) \leq \lambda_1(M(\bfX))$, the definition of $\calL_{\calX, m}$ yields that $\calL_{\calX, m}(g) = 0$ if and only if $\|\Piperp_{V_m(\bfX)} \nabla g(\bfX)\|_F^2 = 0$, which yields the desired equivalence.
\end{proof}

Besides this best-case scenario, we cannot expect in general to have $\calJ_{\calX, m}(g) = 0$ for some $g\in\calG_m$.
A first situation where we can ensure a general result is the bi-Lipschitz case, stated in Proposition \ref{prop:collective bound J by L bi lip} below.

\begin{prpstn}
\label{prop:collective bound J by L bi lip}
    Assume that there exists $c>0$ such that for all $g\in\calG_m$ it holds $c \leq \sigma_{m}(\nabla g(\bfX))^2$ almost surely.
    Then we have
    \[
        \calJ_{\calX, m}(g)
        \leq c^{-1} \calL_{\calX, m}(g).
    \]
\end{prpstn}
\begin{proof}
    This result follows directly from the right inequality from Lemma \ref{lem:bounds on J collective truncated}.
\end{proof}

Note that we lack the reverse bound, as opposed to \cite{nouySurrogatePoincareInequalities2026}.
If we choose to put $\lambda_m(M(\bfX))$ instead of $\lambda_1(M(\bfX))$ in the definition \eqref{equ:def L collective}, then we would straightforwardly obtain a lower bound, but we would lack the upper bound.
In order to obtain both inequalities in Proposition \ref{prop:collective bound J by L bi lip}, or even in the upcoming results, we would need some control on the ratio of eigenvalues $\frac{\lambda_1(M(\bfX))}{\lambda_m(M(\bfX))}$, at least in terms of large deviations.
We leave this for further investigation.

Now if uniform lower bounds are not available for $\sigma_m(\nabla g(\bfX))$, then we shall rely on so-called \emph{small deviations inequalities} or \emph{anti-concentration inequalities}, which consist of upper bounding $\Proba{\sigma_m(\nabla g(\bfX))^2 \leq \alpha}$ for $\alpha >0$, in order to upper bound $\calJ_{\calX, m}$ with $\calL_{\calX, m}$.
Following \cite{nouySurrogatePoincareInequalities2026}, we will assume that the probability measure of $\bfX$ is $s$-concave for $s\in (0,1/d]$, which we define below.

\begin{dfntn}[$s$-concave probability measure]
\label{def:s-concave measures}
  Let $\mu$ be a probability measure on $\bbR^d$ such that $d\mu(\bfx) = \rho(\bfx) d\bfx$.
  For $s\in[-\infty, 1/d]$, $\mu$ is $s$-concave if and only if $\rho$ is supported on a convex set and is
  $\kappa$-concave with $\kappa=s/(1-sd) \in [-1/d, +\infty]$, meaning
  \begin{equation}
    \rho(\lambda \bfx + (1-\lambda) \bfy)
    \geq (\lambda \rho(\bfx)^{\kappa} + (1-\lambda)\rho(\bfy)^{\kappa})^{1/\kappa}
  \end{equation}
  for all $\bfx,\bfy\in\bbR^d$ such that $\rho(\bfx)\rho(\bfy)>0$ and all $\lambda \in [0,1]$.
  The cases $s\in\{-\infty, 0, 1/d\}$ are interpreted by continuity.
\end{dfntn}

An important property of $s$-concave probability measures with $s\in(0,1/d]$ is that they are compactly supported on a convex set.
In particular, a measure is $\frac{1}{d}$-concave if and only if it is uniform.
We refer to \cite{borellConvexSetFunctions1975,borellConvexMeasuresLocally1974} for a deeper study on $s$-concave probability measures.
It is also worth noting that $s$-concave probability measures with $s\in(0,1/d]$ satisfy a Poincar\'e inequality, which is required to obtain \eqref{equ:upper bound approx error nonlinear} for any $u$, although it is not sufficient, as mentioned in Appendix \ref{apx:poincare inequalities for dimension reduction}.

We can now state a small deviation inequality on $\sigma_m(\nabla g(\bfX))^2$ for a polynomial $g$, which is a direct consequence of \cite{nouySurrogatePoincareInequalities2026}, the latter leveraging deviation inequalities from \cite{fradeliziConcentrationInequalitiesSconcave2009}.

\begin{prpstn}
\label{prop:small deviation uniform collective}
    Assume that $\bfX$ is an absolutely continuous random variable on $\bbR^d$ whose distribution is $s$-concave with $s\in(0,1/d]$.
    Assume that $m\geq 2$.
    Let $g:\calX \rightarrow \bbR^m$ be a polynomial with total degree at most $\ell+1 \geq 2$ such that $\Expe{\|\nabla g(\bfX)\|_F^2} \leq m$.
    Then for all $\varepsilon>0$,
    \begin{equation}
        \Proba{
            \sigma_m(\nabla g(\bfX))^2 < q_g \varepsilon
        }
        \leq 2^{5} s^{-1} m^{\frac{1}{4\ell}} \varepsilon^{\frac{1}{2\ell m}},
    \end{equation}
    with $q_g\geq 0$ defined as the median of $\det(\nabla g(\bfX)^T \nabla g(\bfX))$.
\end{prpstn}
\begin{proof}
    The first thing to note is that $\bfx \mapsto \nabla g(\bfx)^T \nabla g(\bfx)$ is a polynomial of total degree at most $2\ell$.
    Then, using \cite[Proposition 3.5]{nouySurrogatePoincareInequalities2026},
    \[
        \Proba{
            \sigma_m(\nabla g(\bfX))^2 < q_g \varepsilon
        }
        \leq 4(1-2^{-s}) s^{-1} 2^{1/4\ell} m^{-1/4\ell}
        \sup_{\bfx\in\calX} \|\nabla g(\bfx)^T \nabla g(\bfx)\|_F^{\frac{m-1}{2\ell m}} \varepsilon^{\frac{1}{2\ell m}},
    \]
    Moreover, we have for all $\bfx\in\calX$,
    \[
        \|\nabla g(\bfx)^T \nabla g(\bfx)\|_F^2 
        = \sum_{i=1}^m \sigma_i(\nabla g(\bfx))^4 
        \leq \left(\sum_{i=1}^m \sigma_i(\nabla g(\bfx))^2 \right)^2
        = \|\nabla g(\bfx)\|_F^4.
    \]
    Also, using \cite[Proposition 3.4]{nouySurrogatePoincareInequalities2026} on $\bfx \mapsto \|\nabla g(\bfx)\|_F^2$, which is also a polynomial of total degree at most $2\ell$, we obtain
    \[
        4^{-2\ell} (1-2^{-s})^{2\ell}
        \sup_{\bfx\in\calX} \|\nabla g(\bfX)\|_F^2
        \leq 
        2 \Expe{\|\nabla g(\bfX)\|_F^2}
        \leq 2 m.
    \]
    Now, by combining the three previous equations and regrouping the exponents we obtain
    \[
        \Proba{
            \sigma_m(\nabla g(\bfX))^2 < q_g \varepsilon
        }
        \leq 2^{4 + \frac{1}{4\ell} - \frac{2}{m} + \frac{m-1}{2\ell m}} s^{-1} (1-2^{-s})^{\frac{1}{m}} m^{\frac{1}{4\ell} (1 - \frac{2}{m})} \varepsilon^{\frac{1}{2\ell m}}.
    \]
    Finally, using $m\geq 2$ and $1-2^{-s} \leq 1$, we obtain the desired result,
    \[
        \Proba{
            \sigma_m(\nabla g(\bfX))^2 < q_g \varepsilon
        }
        \leq 2^{5} s^{-1} m^{\frac{1}{4\ell}} \varepsilon^{\frac{1}{2\ell m}}.
    \]
\end{proof}

Now from the above small deviation inequality, we can upper bound $\calJ_{\calX, m}$ using our surrogate, which we state in Proposition \ref{prop:bound J by L collective} below.

\begin{prpstn}
\label{prop:bound J by L collective}
    Assume that $\bfX$ is an absolutely continuous random variable on $\bbR^d$ whose distribution is $s$-concave with $s\in(0,1/d]$.
    Assume that $m\geq 2$.
    Assume that every $g\in \calG_m$ is a non-constant polynomial with total degree at most $\ell+1 \geq 2$ such that $\Expe{\|\nabla g(\bfX)\|_F^2} \leq m$.
    Assume that $\|\nabla u_Y(\bfX)\|_2 \leq 1$ almost surely.
    Then for all $g\in\calG_m$ and all $p\geq 1$,
    \begin{equation}
        \calJ_{\calX, m}(g)  
        \leq \gamma \nu_{\calG_m, p}^{-\frac{1}{1+2\ell m}} \calL_{\calX, m}(g)^{\frac{1}{1+2\ell m}},
    \end{equation}
    with $\gamma := 2^9 m^{\frac{1}{4\ell}} s^{-1} \min\{s^{-1}, 3p\ell m\}$ and 
    $
        \nu_{\calG_m, p}
        := \inf_{g\in\calG_m} \Expe{\det(\nabla g(\bfX)^T \nabla g(\bfX))^p}^{\frac{1}{p}}.
    $
\end{prpstn}
\begin{proof}
    The proof is similar to the proof of \cite[Proposition 4.5]{nouySurrogatePoincareInequalities2026}.
    Define for all $\alpha>0$ the event $E(\alpha) := (\sigma_m(\nabla g(\bfX))^2 < \alpha)$.
    Then, using that $\|\nabla u_Y(\bfX)\|_2 \leq 1$ almost surely, we obtain
    \[
        \calJ_{\calX, m}(g)
        \leq \Expe{\|\Piperp_{\nabla g(\bfX)} \Pi_{V_m(\bfX)} \nabla u_Y(\bfX)\|_2^2 \mathbbm{1}_{\overline{E(\alpha)}} } 
        + \Proba{E(\alpha)}. 
    \]
    Then, first using the same reasoning as in \eqref{equ: expand Jx}, then using \eqref{equ: bound V proj dg by dg proj V}, then using $\sigma_m(\nabla g(\bfX))^2 \mathbbm{1}_{\overline{E(\alpha)}} \geq \alpha$, and finally using the definition of $\calL_{\calX,m}$ from \eqref{equ:def L collective}, we obtain
    \[
    \begin{aligned}
        \Expe{\|\Piperp_{\nabla g(\bfX)} \Pi_{V_m(\bfX)} \nabla u_Y(\bfX)\|_2^2 \mathbbm{1}_{\overline{E(\alpha)}} } 
        & \leq \Expe{\lambda_1(M(\bfX)) \|\Piperp_{\nabla g(\bfX)} V_m(\bfX)\|_F^2 \mathbbm{1}_{\overline{E(\alpha)}} }
        \\
        & \leq \Expe{\frac{\lambda_1(M(\bfX))}{\sigma_m(\nabla g(\bfX))^2} \| \Piperp_{V_m(\bfX)} \nabla g(\bfX)\|_F^2 \mathbbm{1}_{\overline{E(\alpha)}} }
        \\
        & \leq \alpha^{-1} \Expe{\lambda_1(M(\bfX)) \| \Piperp_{V_m(\bfX)} \nabla g(\bfX)\|_F^2 }
        \\
        & = \alpha^{-1} \calL_{\calX, m}(g).
    \end{aligned}
    \]
    Combining the previous equations with Proposition \ref{prop:small deviation uniform collective} then yields
    \[
        \calJ_{\calX, m}(g)
        \leq \alpha^{-1} \calL_{\calX, m}(g) + \kappa_g \alpha^{\frac{1}{2\ell m}}
        = \kappa_g \left( \kappa_g^{-1} \calL_{\calX, m}(g) \alpha^{-1} + \alpha^{\frac{1}{2\ell m}} \right),
    \]
    with $\kappa_g  := 2^5 s^{-1} m^{\frac{1}{4\ell}} q_g^{-\frac{1}{2\ell m}}$ and $q_g$ as defined in Proposition \ref{prop:small deviation uniform collective}.
    Moreover, from \cite{nouySurrogatePoincareInequalities2026} it holds for any $a\geq 0$ and $b>0$,
    \[
        a^{\frac{b}{1+b}} 
        \leq \inf_{\alpha>0} (a\alpha^{-1} + \alpha^b)
        \leq 2 a^{\frac{b}{1+b}}.
    \]
    Using the above inequality with $a=\kappa_g^{-1} \calL_{\calX, m}(g)$ and $b = 1/2\ell m$, we obtain
    \[
        \calJ_{\calX, m}(g)
        \leq 2 \kappa_g^{1-\frac{1}{1 + 2\ell m}} \calL_{\calX, m}(g)^{\frac{1}{1 + 2\ell m}}
        \leq 2^6 s^{-1} m^{\frac{1}{4\ell}}
        q_g^{-\frac{1}{1+2\ell m}}
        \calL_{\calX, m}(g)^{\frac{1}{1 + 2\ell m}}.
    \]
    Let us now bound $q_g^{-1}$ using moments of $\det(\nabla g(\bfX)^T \nabla g(\bfX))$.
    Using \cite[Proposition 3.4]{nouySurrogatePoincareInequalities2026} on $\bfx \mapsto \det(\nabla g(\bfx)^T \nabla g(\bfx))$ which is a polynomial of total degree at most $2\ell m$, and the fact that $(1-2^{-s})^{-1} \leq 2 s^{-1}$, we obtain
    \[
        q_g^{-1} \leq \Expe{\det(\nabla g(\bfX)^T \nabla g(\bfX))^p}^{-\frac{1}{p}}
        (8 \min\{s^{-1}, 3p\ell m\})^{2\ell m}
        \leq \nu_{\calG_m, p}^{-1}
        \left(8 \min\{s^{-1}, 3p\ell m\}\right)^{2\ell m}.
    \]
    Combining the two previous equations yields the desired result,
    \[
        \calJ_{\calX, m}(g)
        \leq 2^9 m^{\frac{1}{4\ell}} s^{-1} \min\{s^{-1}, 3p\ell m\}
        \nu_{\calG_m, p}^{-\frac{1}{1+2 \ell m}}
        \calL_{\calX, m}(g)^{\frac{1}{1 + 2\ell m}}.
    \]
\end{proof}

It is important to note that the assumption $\Expe{\|\nabla g(\bfX)\|_F^2} \leq m$ is not very restrictive.
For example, it can be satisfied when considering
\begin{equation}
\label{equ:Gm vector space sphere}
    \calG_m = \left\{
        g : \bfx \rightarrow G^T \Phi(\bfx) ~:~
        G \in \bbR^{K\times m}, ~
        G^T \Expe{\nabla \Phi(\bfX)^T \nabla \Phi(\bfX)} G = I_m
    \right\}.
\end{equation}
With this choice of $\calG_m$, the trace identity for the second moment yields $\Expe{\|\nabla g(\bfX)\|_F^2} = m$ for all $g\in\calG_m$.
Note also that one can obtain similar results when $\|\nabla u_Y(\bfX)\|_2 \leq c$ almost surely, for some $c<+\infty$, using the fact that multiplying $u$ by a factor $\alpha$ multiplies both $\calJ_{\calX, m}(g)$ and $\calL_{\calX, m}(g)$ by a factor $\alpha^2$.  
Let us finish this section by pointing out that we encounter the same problem as in \cite{nouySurrogatePoincareInequalities2026}: the exponent in the upper bound in Proposition \ref{prop:bound J by L collective} is $1/(1+2\ell m)$, which scales rather badly with both $m$ and $\ell$ and can be expected to be sharp, as pointed out in \cite{nouySurrogatePoincareInequalities2026}.

\subsection{Minimizing the surrogate}
\label{subsec:minimizing the surrogate collective}

In this section we investigate the problem of minimizing $\calL_{\calX, m}$.
As stated earlier, it is rather straightforward to see that $g \mapsto \calL_{\calX, m}(g)$ is quadratic, which means that we can benefit from various optimization methods from the field of convex optimization. 
In particular, for $\Phi\in\calC^1(\calX, \bbR^K)$ we can express $G \mapsto \calL_{\calX, m}(G^T \Phi)$ as a quadratic form with some positive semidefinite matrix $H_{\calX, m}$ which depends on $u$ and $\Phi$.
This is stated in the following Proposition.

\begin{prpstn}
\label{prop:L collective is quadratic form}
    For any $G\in\bbR^{K\times m}$ it holds 
    \begin{equation}
        \calL_{\calX, m} (G^T \Phi)
        = \Trace{G^T H_{\calX, m} G},
    \end{equation}
    where $H_{\calX, m} := H_{\calX, m}^{(1)} - H_{\calX, m}^{(2)} \in \bbR^{K\times K}$ is a positive semidefinite matrix with
    \begin{equation}
    \begin{aligned}
    \label{equ:def H matrix for L collective}
        H_{\calX,m}^{(1)} &:= \Expe{
            \lambda_1(M(\bfX))
            \nabla \Phi(\bfX)^T \nabla \Phi(\bfX) 
        }, \\
        H_{\calX,m}^{(2)} &:= \Expe{
            \lambda_1(M(\bfX))
            \nabla \Phi(\bfX)^T V_m(\bfX) 
            V_m(\bfX)^T \nabla \Phi(\bfX)
        },
    \end{aligned}
    \end{equation}
    with $M(\bfX)$ and $V_m(\bfX)$ as defined in \eqref{equ:def covariance matrix gradient u} and \eqref{equ:def princial components gradient u}.
\end{prpstn}
\begin{proof}
    First, writing the squared Frobenius norm as a trace and using $(\Piperp_{V_m(\bfX)})^2 = \Piperp_{V_m(\bfX)}$, then 
    switching $\bbE$ with trace, using $\Piperp_{V_m(\bfX)} = I_d - V_m(\bfX) V_m(\bfX)^T$ and using $\nabla g(\bfX) = \nabla \Phi(\bfX) G$, we obtain,
    \begin{align*}
    \begin{aligned}    
        \calL_{\calX, m}(g)
        &= \Expe{\Trace{
            \lambda_1(M(\bfX))
            \nabla g(\bfX)^T 
            \Piperp_{V_m(\bfX)}
            \nabla g(\bfX)
        }}, \\
        &= \Trace{G^T 
            \Expe{
            \lambda_1(M(\bfX))
            \nabla \Phi(\bfX)^T
            (I_d - V_m(\bfX)V_m(\bfX)^T)
            \nabla \Phi(\bfX)
        } G },
    \end{aligned}
    \end{align*}
    which is the desired result.
\end{proof}

As noted in the previous section, the assumption $\Expe{\|\nabla g(\bfX)\|_F^2} \leq m$ from Proposition \ref{prop:bound J by L collective} can be satisfied by considering $\calG_m$ of the form
\begin{equation}
    \calG_m = \left\{
        g : \bfx \rightarrow G^T \Phi(\bfx) ~:~
        G \in \bbR^{K\times m}, ~
        G^T R G = I_m
    \right\}
\end{equation}
with $R := \Expe{\nabla \Phi(\bfX)^T \nabla \Phi(\bfX)}\in\bbR^{K\times K}$ a symmetric positive definite matrix.
Note that as pointed out in \cite{bigoniNonlinearDimensionReduction2022}, the orthogonality condition $G^T R G = I_m$ has no impact on the minimization of $\calJ_{\calX}$ or its truncated version $\calJ_{\calX, m}$, because $\Pi_{\nabla g(\bfX)}$ is invariant to invertible transformations of $g$.
In this context, minimizing $\calL_{\calX, m}$ over $\calG_m$ is equivalent to finding the minimal generalized eigenpair of the pencil $(H_{\calX, m}, R)$, as stated in Proposition \ref{prop:L collective as generalized ev}.

\begin{prpstn}
\label{prop:L collective as generalized ev}
    Let $\calG_m$ be as in \eqref{equ:Gm vector space sphere}.
    The minimizers of $\calL_{\calX, m}$ over $\calG_m$ are the functions of the form $g^*(\bfx) = (G^*)^T \Phi(\bfx)$, where $G^* \in \bbR^{K\times m}$ is a solution to the generalized eigenvalue problem
    \begin{equation}
    \label{equ:L collective as generalized ev}
        \min_{\substack{G \in \bbR^{K\times m} \\G^T R G = I_m}} 
        \Trace{G^T H_{\calX,m} G},
    \end{equation}
    with $H_{\calX,m}$ defined in \eqref{equ:def H matrix for L collective}.
\end{prpstn}

We end this section by discussing the major computational problem with $\calL_{\calX,m}$.
Indeed, while $\calJ_{\calX}(g)$ can be estimated by classical Monte-Carlo methods by independently sampling $(\bfx^{(i)}, y^{(i)})_{1\leq i\leq n_s}$ from $\mu_{\calX} \otimes \mu_{\calY}$, this is not the case for $\calL_{\calX,m}(g)$ as it requires estimating $\lambda_1(M(\bfx^{(i)}))$ and $V_m(\bfx^{(i)})$ for all samples.
One way to do so is to use a tensorized sample $(\bfx^{(i)}, y^{(j)})_{1\leq i \leq n_{\calX}, 1\leq j \leq n_{\calY}}$, of size $n_s = n_{\calX}n_{\calY}$.


\section{Two groups setting}
\label{sec:two variables approach}

In this section, we consider $\bfX$ with measure $\mu$ over $\calX \subset \bbR^d$.
We fix a multi-index $\alpha \subset \{1, \cdots, d\}$, and we assume that $\bfX_{\alpha} := (X_i)_{i\in\alpha}$ and $\bfX_{\alpha^c}:=(X_{i})_{i\in\alpha^c}$ are independent, meaning that $\mu = \mu_{\alpha} \otimes \mu_{\alpha^c}$ with support $\calX_{\alpha} \times \calX_{\alpha^c}$.
For any strictly positive integers $n_{\alpha}$ and $n_{\alpha^c}$, and any functions $h^{\alpha} : \calX_{\alpha} \rightarrow \bbR^{n_{\alpha}}$ and $h^{\alpha^c} : \calX_{\alpha^c} \rightarrow \bbR^{n_{\alpha^c}}$, we identify the tuple $(h^{\alpha}, h^{\alpha^c})$ with the function $\bfx \mapsto (h^{\alpha}(\bfx_{\alpha}), h^{\alpha^c}(\bfx_{\alpha^c})) \in \bbR^{n_{\alpha} + n_{\alpha^c}}$.

For some fixed $\bfm = (m_{\alpha}, m_{\alpha^c}) \in \bbN\times \bbN$ and fixed classes of functions $\calF_{\bfm}$ and $ \calG_{m_{\alpha}}^{\alpha}$, and $\calG_{m_{\alpha^c}}^{\alpha^c}$, we then consider an approximation of the form $\bfx \mapsto f \circ g(\bfx)$, with some regression function $f : \bbR^{m_{\alpha}} \times \bbR^{m_{\alpha^c}} \rightarrow \bbR$ from $\calF_{\bfm}$ and some separated feature map $(g^{\alpha}, g^{\alpha^c}) \equiv g : \calX \rightarrow \bbR^{m_{\alpha}} \times \bbR^{m_{\alpha^c}}$ from $\calG_{\bfm} \equiv \calG^{\alpha}_{m_{\alpha}} \times \calG^{\alpha^c}_{m_{\alpha^c}}$ such that 
\[
  g: \bfx \mapsto (g^{\alpha}(\bfx_{\alpha}), g^{\alpha^c}(\bfx_{\alpha^c})),
\]
with $g^{\alpha} : \calX_{\alpha} \rightarrow \bbR^{m_{\alpha}}$ from $\calG_{m_{\alpha}}^{\alpha}$ and $g^{\alpha^c} : \calX_{\alpha^c} \rightarrow \bbR^{m_{\alpha^c}}$ from $\calG_{m_{\alpha^c}}^{\alpha^c}$.
We are then considering
\begin{equation}
\label{equ:approximation error bivariate}
  \inf_{g \in \calG_{\bfm}} 
  \inf_{f\in\calF_{\bfm}}
  \Expe{|u(\bfX) - f(g^{\alpha}(\bfX_{\alpha}), g^{\alpha^c}(\bfX_{\alpha^c}))|^2}.
\end{equation}

In this section we discuss different approaches for solving or approximating \eqref{equ:approximation error bivariate} depending on choices for $\calF_{\bfm}$.
First in \Cref{subsec:bilinear regression function} we discuss bilinear regression functions, which is related to classical singular value decomposition.
Then in \Cref{subsec:unconstrained regression function bi} we discuss unconstrained regression functions, assuming only measurability, which corresponds to a more general dimension reduction framework.

\subsection{Bilinear regression function}
\label{subsec:bilinear regression function}

In this section we discuss the case where $\calF_{\bfm} = \calF_{\bfm}^{bi}$ contains only bilinear functions, in the sense that $f(\bfz^{\alpha},\cdot)$ and $f(\cdot, \bfz^{\alpha^c})$ are linear for any $(\bfz^{\alpha},\bfz^{\alpha^c}) \in \bbR^{m_{\alpha}} \times \bbR^{m_{\alpha^c}}$.
In other words we identify $\calF_{\bfm}^{bi}$ with $\bbR^{m_{\alpha}\times m_{\alpha^c}}$, and we want to minimize over $\calG_{\bfm}$ the function
\begin{align}
\begin{aligned}
\label{equ:approximation error bilinear}
  \calE_{\alpha}^{bi} : g
  \mapsto  & \inf_{A\in\bbR^{m_{\alpha}\times m_{\alpha^c}}}
  \Expe{|u(\bfX) - g^{\alpha}(\bfX_{\alpha})^T A g^{\alpha^c}(\bfX_{\alpha^c})|^2}
  \\
  & = \inf_{A\in\bbR^{m_{\alpha}\times m_{\alpha^c}}}
  \|u - \sum_{i,j} A_{ij}
    g_i^\alpha \otimes g_j^{\alpha^c}
  \|_{L^2(\mathcal X,\mu)}^2.
\end{aligned}
\end{align}
For fixed $g \in \calG_{\bfm}$ with $g\equiv (g^{\alpha}, g^{\alpha^c})$, the optimal $A\in\bbR^{m_{\alpha}\times m_{\alpha^c}}$ is given via the orthogonal projection of $u$ onto the subspace 
\[
  \spanv{g^{\alpha}_i \otimes g^{\alpha^c}_j : 1\leq i\leq m_{\alpha}, 1\leq j\leq m_{\alpha^c}}
  =
  \spanv{g^{\alpha}_i}_{1\leq i\leq m_{\alpha}} 
  \otimes 
  \spanv{g^{\alpha^c}_j}_{1\leq j\leq m_{\alpha^c}}
\]
with $A_{ij} = \innerp{g^{\alpha}_i \otimes g^{\alpha^c}_j}{u}$ when $(g_i^{\alpha} \otimes g^{\alpha^c}_j)_{1\leq i\leq m_{\alpha}, 1\leq j\leq m_{\alpha^c}}$ are orthonormal in $L^2(\calX, \mu)$.
Note that \eqref{equ:approximation error bilinear} is actually invariant to any invertible linear transformation of elements of $\calG_{m_{\alpha}}^{\alpha}$ and $\calG_{m_{\alpha^c}}^{\alpha^c}$, meaning that it only depends on $U_{\alpha}=\spanv{g^{\alpha}_i}_{1\leq i\leq m_{\alpha}}$ and $U_{\alpha^c}=\spanv{g^{\alpha^c}_j}_{1\leq j\leq m_{\alpha^c}}$.

Now assume that $\calG^{\alpha}_{m_{\alpha}}$ and $\calG^{\alpha^c}_{m_{\alpha^c}}$ are vector spaces  such that the components of $g^{\alpha}$ and $g^{\alpha^c}$ lie respectively in some fixed vector spaces $V_{\alpha} \subset L^2(\calX_{\alpha}, \mu_{\alpha})$ and $V_{\alpha^c} \subset L^2(\calX_{\alpha^c}, \mu_{\alpha^c})$.
In this case, the optimal $g\in\calG_{\bfm}$ is given via the singular value decomposition of $\calP_{V_{\alpha} \otimes V_{\alpha^c}} u$, see for example \cite[Section 4.4.3]{hackbuschTensorSpacesNumerical2019}.
This decomposition is written as
\begin{equation}
\label{equ:def svd}
  (\calP_{V_{\alpha} \otimes V_{\alpha^c}} u)(\bfx) 
  = \sum_{k=1}^{\min(\dim V_{\alpha}, ~ \dim V_{\alpha^c})} 
  \sigma_k^{\alpha} v^{\alpha}_k(\bfx_{\alpha}) v^{\alpha^c}_k(\bfx_{\alpha^c}),
\end{equation}
where $(v^{\alpha}_i)_{1 \leq i \leq \dim V_{\alpha}}$ and $(v^{\alpha^c}_j)_{1 \leq j \leq \dim V_{\alpha^c}}$ are singular vectors, which form orthonormal bases of $V_{\alpha}$ and $V_{\alpha^c}$ respectively, with associated singular values $\sigma_1^{\alpha} \geq \sigma_2^{\alpha} \geq \cdots$.
Then the optimal $g\in\calG_{\bfm}$ is obtained by truncating the above sum, keeping only the first $\min(m_{\alpha},m_{\alpha^c})$ terms, which reads
\[
  \hat u(\bfx) 
  = \sum_{k=1}^{\min(m_{\alpha},m_{\alpha^c})} 
  \sigma^{\alpha}_k 
  v^{\alpha}_k(\bfx_{\alpha}) 
  v^{\alpha^c}_k(\bfx_{\alpha^c}).
\]
In particular, there are only $\min(m_{\alpha}, m_{\alpha^c})$ terms in the sum.
Thus, it is equivalent to consider $m_{\alpha} = m_{\alpha^c}$.
Finally, a minimizer of \eqref{equ:approximation error bilinear} is given by $g^{\alpha} = (v^{\alpha}_i)_{1\leq i\leq m_{\alpha}}$, $g^{\alpha^c} = (v^{\alpha^c}_i)_{1\leq i\leq m_{\alpha}}$ and $A = \mathrm{diag}((\sigma^{\alpha}_i)_{1\leq i\leq m_{\alpha}})$.
Also, if the singular values are all distinct, then $\spanv{g^{\alpha}}$ and $\spanv{g^{\alpha^c}}$ are unique.
The associated approximation error \eqref{equ:approximation error bivariate} is given by 
\[
  \min_{g\in\calG_{\bfm}} 
  \calE_{\alpha}^{bi}(g)
  =
  \|u - \calP_{V_{\alpha} \otimes V_{\alpha^c}}u\|^2_{L^2} 
  + \sum_{k=m_{\alpha} + 1}^{\min(\dim V_{\alpha}, ~ \dim V_{\alpha^c})} (\sigma_k^{\alpha})^2.
\]

Let us emphasize the fact that, due to the SVD truncation property, the resulting number of features is the same for both $\bfX_{\alpha}$ and $\bfX_{\alpha^c}$.
This is an interesting feature of SVD-based approximation, as low dimensionality with respect to $\calX_{\alpha}$ implies low dimensionality with respect to $\calX_{\alpha^c}$, and vice versa. 
This is also interesting for practical algorithms as the singular vectors in $V_{\alpha}$ can be estimated independently of those in $V_{\alpha^c}$.
For example, when $\dim \calX_{\alpha}$ is much smaller than $\dim \calX_{\alpha^c}$, sampling-based estimation is easier for $v_k^{\alpha}$ than for $v_k^{\alpha^c}$.

We end this section by noting that this bilinear framework will also be relevant in the multilinear framework discussed in \Cref{sec:multiple variables approach}, especially the optimality of SVD.

\subsection{Unconstrained regression function}
\label{subsec:unconstrained regression function bi}

In this section we discuss the case where there is no restriction besides measurability on $\calF_{\bfm}$, meaning that $\calF_{\bfm} \equiv \calF_m := \{f : \bbR^{m} \rightarrow \bbR ~ \text{measurable}\}$ with $m := m_{\alpha} + m_{\alpha^c}$.
We then want to minimize over $\calG_{\bfm} \equiv \calG^{\alpha}_{m_{\alpha}} \times \calG^{\alpha^c}_{m_{\alpha^c}}$ the function $\calE$ defined for any $g \equiv (g^{\alpha}, g^{\alpha^c})$ by 
\begin{equation}
\label{equ:def E}
  \calE(g) := 
  \inf_{f\in \calF_m}
  \Expe{|u(\bfX) - f(g^{\alpha}(\bfX_{\alpha}), g^{\alpha^c}(\bfX_{\alpha^c}))|^2}.
\end{equation}
The function $f\in\calF_m$ satisfying the above infimum is given via an orthogonal projection onto some subspace of $L^2(\calX, \mu)$, the subspace of $g$-measurable functions
\begin{equation}
\label{equ:def Sigma g}
  \Sigma(g) 
  := L^2(\calX, \sigma(g(\bfX)), \mu) 
  = \{\bfx \mapsto f(g^{\alpha}(\bfx_{\alpha}), g^{\alpha^c}(\bfx_{\alpha^c})): ~ f\in\calF_m\}
  \cap  L^2(\calX, \mu).
\end{equation}
The function $f$ associated to the projection of $u$ onto $\Sigma(g)$ is given via the conditional expectation $f(\bfz^{\alpha}, \bfz^{\alpha^c}) = \Expe{u(\bfX) | g(\bfX) = (\bfz^{\alpha}, \bfz^{\alpha^c})}$.
Moreover, since $\mu = \mu_{\alpha} \otimes \mu_{\alpha^c}$, the subspace $\Sigma(g)$ is a tensor product, $\Sigma(g) = \Sigma_{\alpha}(g^{\alpha}) \otimes \Sigma_{\alpha^c}(g^{\alpha^c})$, where for $\beta \subset \{1, \cdots, d\}$,
\begin{equation}
\label{equ:def Sigma beta g}
  \Sigma_{\beta}(g^{\beta}) 
  := L^2(\calX_{\beta}, \sigma(g^{\beta}(\bfX_{\beta})), \mu_{\beta}) 
  = \{h\circ g^{\beta} : h : \bbR^{m_{\beta}} \rightarrow \bbR \text{ measurable}\}
  \cap  L^2(\calX_{\beta}, \mu_{\beta}).
\end{equation}

There are several differences compared to the bilinear case.
A first difference is that $\Sigma(g)$ is an infinite dimensional space, contrary to $U_{\alpha} \otimes U_{\alpha^c} = \spanv{g^{\alpha}_i \otimes g^{\alpha^c}_j}_{i,j}$.
Hence for building $f$ in practice, we approximate $\Sigma(g)$ by a finite dimensional space.
A second difference is that if $g^{\alpha}$ reproduces the identity, meaning that $Rg^{\alpha}(\bfX_{\alpha}) = \bfX_{\alpha}$ for some matrix $R \in \bbR^{|\alpha| \times m_{\alpha}}$, then $\Sigma_{\alpha}(g^{\alpha}) = \Sigma_{\alpha}(id^{\alpha}) = L^2(\calX_{\alpha}, \mu_{\alpha})$.
The same holds for $g^{\alpha^c}$.
This means that taking $m_{\alpha} \geq  |\alpha|$ or $m_{\alpha^c} \geq |\alpha^c|$ is somewhat useless in this setting.
A third difference is that, even with strong assumptions on $\calG^{\alpha}_{m_{\alpha}}$ and $\calG^{\alpha^c}_{m_{\alpha^c}}$, minimization of $\calE$ over $\calG_{\bfm} \equiv \calG^{\alpha}_{m_{\alpha}} \times \calG^{\alpha^c}_{m_{\alpha^c}}$ is not related to a classical approximation problem, such as SVD.
This is a crucial difference as optimality in the two groups setting can be leveraged to obtain near-optimality in the multiple groups setting, as discussed in \Cref{sec:multiple variables approach}.

Hence, as in the one variable framework, we can only consider heuristics or upper bounds on $\calE$ to obtain suboptimal $g$.
For example, when considering Poincar\'e inequality-based methods, the product structure of $\calG_{\bfm} \equiv \calG^{\alpha}_{m_{\alpha}} \times \calG^{\alpha^c}_{m_{\alpha^c}}$ transfers naturally to $\calJ$, as stated in Proposition \ref{prop:poincare bound decomposition bi} below.

\begin{prpstn}
\label{prop:poincare bound decomposition bi}
  For any $g \in \calG_{\bfm}$ with $g \equiv (g^{\alpha}, g^{\alpha^c})$, it holds
  \begin{equation}
    \calJ(g) 
    = \calJ((g^{\alpha}, id^{\alpha^c}))
    + \calJ((id^{\alpha}, g^{\alpha^c}))
    = \calJ_{\calX_{\alpha}}(g^{\alpha})
    + \calJ_{\calX_{\alpha^c}}(g^{\alpha^c}),
  \end{equation}
  with $\calJ_{\calX_{\alpha}}$ and $\calJ_{\calX_{\alpha^c}}$ as defined in \eqref{equ:def of J collective}.
\end{prpstn}
\begin{proof}
  We refer to the more general proof of Proposition \ref{prop:poincare bound decomposition multi}.
\end{proof}

A consequence of Proposition \ref{prop:poincare bound decomposition bi} is that minimizing $\calJ$ over $\calG_{\bfm} \equiv \calG_{m_{\alpha}}^{\alpha} \times \calG_{m_{\alpha^c}}^{\alpha^c}$ is equivalent to minimizing $\calJ_{\calX_{\alpha}}$ and $\calJ_{\calX_{\alpha^c}}$ over $\calG_{m_{\alpha}}^{\alpha}$ and $\calG_{m_{\alpha^c}}^{\alpha^c}$ respectively.
In other words, it is equivalent to solving two instances of the Poincar\'e inequality-based collective setting described in \Cref{sec:collective dimension reduction}.
As a result, one may consider leveraging the surrogates $\calL_{\calX_{\alpha}, m_{\alpha}}$ and $\calL_{\calX_{\alpha^c}, m_{\alpha^c}}$ from \Cref{sec:collective dimension reduction}.
Note also that the same tensorized sample can be used for both surrogates.


\section{Multiple groups setting}
\label{sec:multiple variables approach}

In this section we fix $S$ a partition of $D :=\{1, \cdots, d\}$ of size $N>1$, meaning that $S:=\{\alpha_1, \cdots, \alpha_N\}$ where $\bigcup_{\alpha \in S} \alpha = \{1, \cdots, d\}$ and the union is disjoint.
We assume that $(\bfX_{\alpha})_{\alpha \in S}$ are independent random vectors, meaning that $\mu = \otimes_{\alpha \in S}\mu_{\alpha}$.
For any strictly positive integers $(n_{\alpha})_{\alpha \in S}$ and any functions $h^{\alpha} : \calX_{\alpha} \rightarrow \bbR^{n_{\alpha}}$, we identify the tuple $(h^{\alpha})_{\alpha \in S}$ with the function $\bfx \mapsto (h^{\alpha_1}(\bfx_{\alpha_1}), \cdots, h^{\alpha_{N}}(\bfx_{\alpha_{N}})) \in \bbR^{n_{\alpha_1} + \cdots + n_{\alpha_{N}}}$.

For some fixed $\bfm = (m_{\alpha})_{\alpha \in S}\in\bbN^N$ and fixed classes of functions $\calF_{\bfm}$ and $(\calG_{m_{\alpha}}^{\alpha})_{\alpha \in S}$, we then discuss an approximation of the form $\bfx \mapsto f \circ g(\bfx)$, with some regression function $f : \times_{\alpha\in S}\bbR^{m_{\alpha}} \rightarrow \bbR$ from $\calF_{\bfm}$ and some separated feature map $(g^{\alpha})_{\alpha \in S} \equiv g : \calX \rightarrow \times_{\alpha \in S} \bbR^{m_{\alpha}}$ from $\calG_{\bfm} \equiv \times_{\alpha \in S} \calG^{\alpha}_{m_{\alpha}}$, such that 
\[
  g(\bfx) = (g^{\alpha_1}(\bfx_{\alpha_1}), \cdots, g^{\alpha_N}(\bfx_{\alpha_N})),
\]
with $g^{\alpha} : \calX_{\alpha} \rightarrow \bbR^{m_{\alpha}}$ from $\calG^{\alpha}_{m_{\alpha}}$ for all $\alpha \in S$.
We are then considering
\begin{equation}
\label{equ:approximation error multi}
  \inf_{g \in \calG_{\bfm}} 
  \inf_{f\in\calF_{\bfm}}
  \Expe{|u(\bfX) - f(g^{\alpha_1}(\bfX_{\alpha_1}), \cdots, g^{\alpha_N}(\bfX_{\alpha_N}))|^2}.
\end{equation}
In this section we discuss different approaches for tackling \eqref{equ:approximation error multi} depending on choices for $\calF_{\bfm}$.
In \Cref{subsec:multilinear regression function} we discuss multilinear regression functions, which corresponds to tensor-based approximation in Tucker format.
Then in \Cref{subsec:unconstrained regression function multi} we discuss unconstrained measurable regression functions, assuming only measurability, which corresponds to a more general dimension reduction framework.

\subsection{Multilinear regression function}
\label{subsec:multilinear regression function}

In this section we discuss the case where $\calF_{\bfm} = \calF_{\bfm}^{mul}$ contains only multilinear functions, in the sense that for all $\alpha \in S$ and all $(\bfz^{\beta})_{\beta \in S \setminus \{\alpha\}}$, the function $f(\cdot, (\bfz^{\beta})_{\beta \in S \setminus \{\alpha\}})$ is linear.
In other words $\calF_{\bfm}^{mul} \equiv \bbR^{m_{\alpha_1} \times \cdots \times m_{\alpha_N}}$ is a set of tensors of order $N$, and the composition with $g$ can be written as
\[
  f\circ g(\bfx)
  = \sum_{i_1, \cdots, i_N}  T_{i_1, \cdots, i_N} g_{i_1}^{\alpha_1}(\bfx_{\alpha_1}) \cdots g_{i_N}^{\alpha_N}(\bfx_{\alpha_N}).
\]
We then want to minimize over $\calG_{\bfm}$ the function
\begin{align}
\begin{aligned}
\label{equ:approximation error multilinear}
  \calE_{S}^{mul} : g
  \mapsto
  & \inf_{T \in \calF_{\bfm}^{mul}}
  \Expe{|u(\bfX) - T((g^{\alpha}(\bfX_{\alpha}))_{\alpha \in S})|^2}
  \\
  & = \inf_{T\in\bbR^{m_{\alpha_1}\times \cdots \times m_{\alpha_N}}}
  \|u - \sum_{i_1, \cdots, i_N}  T_{i_1, \cdots, i_N} g_{i_1}^{\alpha_1}\otimes \cdots \otimes g_{i_N}^{\alpha_N} \|_{L^2(\mathcal X,\mu)}^2.
\end{aligned}
\end{align}
For fixed $g\in \calG_{\bfm}$ with $g \equiv (g^{\alpha})_{\alpha \in S}$, the optimal tensor $T$ is given via the orthogonal projection of $u$ onto the subspace 
\[
  \bigotimes_{\alpha \in S}
  \spanv{g^{\alpha}_i}_{1\leq i\leq m_{\alpha}}
\]
with $T_{(i_{\alpha})_{\alpha \in S}} = \innerp{\otimes_{\alpha \in S} g_{i_{\alpha}}^{\alpha} }{u}$ when the $(\otimes_{\alpha \in S}g^{\alpha}_{i_{\alpha}})$ are orthonormal in $L^2(\calX, \mu)$.
Similarly to the bilinear case, we can again note that for each $\alpha \in S$, \eqref{equ:approximation error multilinear} is actually invariant to any invertible linear transformation on elements of $\calG_{m_{\alpha}}^{\alpha}$, meaning that it only depends on $U_{\alpha}=\spanv{g^{\alpha}_i}_{1\leq i\leq m_{\alpha}}$.

Now assume that for every $\alpha \in S$, $\calG_{m_{\alpha}}^{\alpha}$ is a vector space such that the components of $g^{\alpha}$ lie in some fixed vector spaces $V_{\alpha} \subset L^2(\calX_{\alpha}, \mu_{\alpha})$.
This setting actually corresponds to the so-called \emph{tensor subspace} (or \emph{Tucker}) format \cite[Chapter 10]{hackbuschTensorSpacesNumerical2019}, and comes  with multiple optimization methods for minimizing $\calE_S^{mul}$ over $\calG_{\bfm}$.
We will focus on the so-called higher-order singular value decomposition (HOSVD), which is defined for all $\alpha \in S$ by $g^{\alpha}_{\mathrm{HOSVD}} = (v^{\alpha}_1, \cdots, v^{\alpha}_{m_{\alpha}})$, with $v^{\alpha}_k$ as defined in \eqref{equ:def svd} 
with $V_{\alpha^c} = \bigotimes_{\beta \in S \setminus \{\alpha\}} V_{\beta}$, which is optimal with respect to $\calE^{bi}_{\alpha}$ defined in \eqref{equ:approximation error bilinear}.
Then, with $g_{\mathrm{HOSVD}} \equiv (g^{\alpha}_{\mathrm{HOSVD}})_{\alpha \in S}$, \cite[Theorem 10.2]{hackbuschTensorSpacesNumerical2019} states that
\[
  \inf_{T\in \calF_{\bfm}^{mul}} 
  \| \calP_{\otimes_{\alpha \in S} V_{\alpha}} u - T \circ g_{\mathrm{HOSVD}}\|_{L^2}^2
  \leq 
  N
  \inf_{g \in \calG_{\bfm}}
  \inf_{T \in \calF_{\bfm}^{mul}}
  \| \calP_{\otimes_{\alpha \in S} V_{\alpha}} u - T\circ g\|_{L^2}^2,
\]
in other words that $g_{\mathrm{HOSVD}}$ is near-optimal.
Moreover, since for all $T\in\calF_{\bfm}^{mul}$ and all $g\in\calG_{\bfm}$ we have $T\circ g \in \bigotimes_{\alpha \in S} V_{\alpha}$, we have that
\[
  \calE_S^{mul}(g) - 
  \|u - \calP_{\otimes_{\alpha \in S} V_{\alpha}} u\|_{L^2}^2 
  = \inf_{T\in\calF^{mul}_{\bfm}} 
  \|\calP_{\otimes_{\alpha \in S} V_{\alpha}} u - T \circ g\|_{L^2}^2.
\]
As a result, by first applying the above equality to $g_{\mathrm{HOSVD}}$, then applying the quasi-optimality result of the HOSVD, then taking the infimum over $g$ in the above equality, we obtain
\begin{align}
\begin{aligned}
\label{equ:near optimal multilinear}
  \calE_S^{mul}(g_{\mathrm{HOSVD}}) 
  - \|u - \calP_{\otimes_{\alpha \in S} V_{\alpha}} u\|^2_{L^2}
  & \leq N
    \inf_{g \in \calG_{\bfm}}
    \inf_{T \in \calF_{\bfm}^{mul}}
    \| \calP_{\otimes_{\alpha \in S} V_{\alpha}} u - T\circ g\|_{L^2}^2
  \\
  & \leq
  N \left( 
    \inf_{g\in \calG_{\bfm}}\calE_S^{mul}(g) 
    -\|u - \calP_{\otimes_{\alpha \in S} V_{\alpha}} u\|^2_{L^2}
    \right).
\end{aligned}
\end{align}

\subsection{Unconstrained regression function}
\label{subsec:unconstrained regression function multi}

In this section we discuss the case where there is no restriction besides measurability on $\calF_{\bfm} = \calF_m$, meaning that $\calF_m = \{f : \bbR^m \rightarrow \bbR ~ \text{measurable}\}$ with $m := \sum_{\alpha \in S} m_{\alpha}$.
We then want to minimize over $\calG_{\bfm} \equiv \times_{\alpha \in S} \calG^{\alpha}_{m_{\alpha}}$ the function $\calE$ defined for any $g\equiv (g^{\alpha})_{\alpha \in S}$ by
\begin{equation}
\label{equ:approximation error unconstrained multi}
  \calE(g) =
  \inf_{f\in\calF_m}
  \Expe{|u(\bfX) - f(g^{\alpha_1}(\bfX_{\alpha_1}), \cdots, g^{\alpha_N}(\bfX_{\alpha_N}))|^2}.
\end{equation}
For fixed $g\in \calG_{\bfm}$ with $g \equiv(g^{\alpha})_{\alpha \in S}$, the optimal $f\in\calF_m$ is again given via an orthogonal projection onto $\Sigma(g)$, given via the conditional expectation $f(\bfz) = \Expe{u(\bfX) | g(\bfX) = \bfz}$.
Moreover, since $\mu = \otimes_{\alpha \in S} \mu_{\alpha}$, the subspace $\Sigma(g)$ is again a tensor product, $\Sigma(g) = \otimes_{\alpha \in S} \Sigma_{\alpha}(g^{\alpha})$.

The fact that $\calE(g)$ is a projection error onto a tensor product space allows us to make a link with the two groups setting from \Cref{subsec:unconstrained regression function bi}, similarly to HOSVD.
In particular, the optimization on $\calG_{\bfm}$ is nearly equivalent to $N$ separate optimization problems on $\calG_{m_{\alpha}}^{\alpha}$ for $\alpha \in S$.
This is stated in Proposition \ref{prop:approximation error decomposition multi} below.
For any $\beta\subset\{1, \cdots, N\}$ we write $id^{\beta}$ the identity function on $\bbR^{|\beta|}$.

\begin{prpstn}
\label{prop:approximation error decomposition multi}
  Assume that $\mu = \otimes_{\alpha \in S} \mu_{\alpha}$. 
  Then for all $g \in \calG_{\bfm}$ with $g\equiv (g^{\alpha})_{\alpha\in S}$, it holds 
  \begin{equation}
  \label{equ:approximation error decomposition multi}
    \calE(g) 
    \leq 
    \sum_{\alpha \in S} 
    \calE((g^{\alpha}, id^{\alpha^c}))
    =  \sum_{\alpha \in S} 
    \calE_{\calX_{\alpha}}(g^{\alpha})
    \leq
    N \calE(g),
  \end{equation}
  with $\calE_{\calX_{\alpha}}$ as defined in \eqref{equ:def reconstruction error collective}.
\end{prpstn}
\begin{proof}
  Firstly, for any $\alpha \in S$ we have $\Sigma(g) \subset \Sigma((g^{\alpha}, id^{\alpha^c}))$, thus $\calE((g^{\alpha}, id^{\alpha^c})) \leq \calE(g)$, where $\calE((g^{\alpha}, id^{\alpha^c})) = \calE_{\calX_{\alpha}}(g^{\alpha})$.
  Summing those inequalities for all $\alpha \in S$ yields the desired right inequality in \eqref{equ:approximation error decomposition multi}.
  Secondly, the product structure of $\mu$ implies that $\calP_{\Sigma(g)} = \prod_{\alpha \in S} \calP_{\Sigma((g^{\alpha}, id^{\alpha^c}))}$, where the projectors in the right-hand side commute.
  Now from \cite[Lemma 4.145]{hackbuschTensorSpacesNumerical2019} it holds that
  \[
    \calE(g)  
    =\|(I - \prod_{\alpha \in S} \calP_{\Sigma((g^{\alpha}, id^{\alpha^c}))} ) u\|_{L ^2}^2
    \leq
    \sum_{\alpha \in S}
    \|\calP^{\perp}_{\Sigma((g^{\alpha}, id^{\alpha^c}))}u\|_{L ^2}^2
     = \sum_{\alpha \in S} \calE((g^{\alpha}, id^{\alpha^c})).
  \]
  This yields the desired left inequality in \eqref{equ:approximation error decomposition multi}, which concludes the proof.
\end{proof}

A direct consequence of Proposition \ref{prop:approximation error decomposition multi} is that minimizers of $\calE_{\calX_{\alpha}}$ over $\calG^{\alpha}_{m_{\alpha}}$ for $\alpha \in S$, if they actually exist, are near-optimal when minimizing $\calE$ over $\calG_{\bfm}$.
This is stated in Corollary \ref{coro:approximation error quasi optimal multi}, and is similar to the near optimality result \eqref{equ:near optimal multilinear}.

\begin{crllr}
\label{coro:approximation error quasi optimal multi}
  Assume that $\mu = \otimes_{\alpha \in S} \mu_{\alpha}$, and that for all $\alpha \in S$ there exists a minimizer $g_*^{\alpha}$ of $\calE((\cdot, id^{\alpha^c}))$ over $\calG_{m_{\alpha}}^{\alpha}$.
  Then for $g_* \equiv (g^{\alpha}_*)_{\alpha \in S}$ it holds 
  \[
    \calE(g_*) 
    \leq N \inf_{g \in \calG_{\bfm}} \calE(g).
  \]
\end{crllr}
\begin{proof}
  Let $g \in \calG_{\bfm}$ with $g \equiv (g^{\alpha})_{\alpha \in S}$.
  Using the left inequality from \eqref{equ:approximation error decomposition multi}, then using the definition of $(g^{\alpha}_*)_{\alpha \in S}$ and the right inequality from \eqref{equ:approximation error decomposition multi}, we obtain
  \begin{equation}
    \calE(g_*) 
    \leq
    \sum_{\alpha \in S} 
    \calE((g^{\alpha}_*, id^{\alpha^c}))
    \leq
    \sum_{\alpha \in S} 
    \calE((g^{\alpha}, id^{\alpha^c}))
    \leq  
    N \calE(g) .
  \end{equation}
\end{proof}

Unfortunately, while the HOSVD from the multilinear case in \Cref{subsec:multilinear regression function} leverages the fact that a minimizer of $\calE_{\alpha}^{bi}$ is given by the SVD, here the minimization of $\calE((\cdot, id^{\alpha^c})) = \calE_{\calX_{\alpha}}(\cdot)$ remains a challenge.
Hence, we can only consider heuristics or upper bounds on the latter, as investigated in \Cref{sec:collective dimension reduction}.
For example, when considering Poincaré inequality-based methods, as in \Cref{sec:two variables approach}, the product structure of $\calG_{\bfm} \equiv \times_{\alpha \in S} \calG^{\alpha}_{m_{\alpha}}$ transfers naturally to $\calJ$ by its definition, as stated in Proposition \ref{prop:poincare bound decomposition multi} which generalizes Proposition \ref{prop:poincare bound decomposition bi} for the two groups setting.
For any $\beta\subset\{1, \cdots, N\}$ we write $\nabla_{\beta} u(\bfx) = (\frac{\partial u}{\partial x_i}(\bfx))_{i \in \beta} \in \bbR^{|\beta|}$.

\begin{prpstn}
\label{prop:poincare bound decomposition multi}
  For any $g = (g^{\alpha})_{\alpha \in S} \in \calG_{\bfm}$, it holds
  \begin{equation}
    \calJ(g) 
    = \sum_{\alpha \in S}
    \Expe{\|\Piperp_{\nabla g^{\alpha}(\bfX_{\alpha})} \nabla_{\alpha} u(\bfX)\|_2^2}
    = \sum_{\alpha \in S}
    \calJ((g^{\alpha}, id^{\alpha^c}))
    = \sum_{\alpha \in S}
    \calJ_{\calX_{\alpha}}(g^{\alpha}),
  \end{equation}
  with $\calJ_{\calX_{\alpha}}$ as defined in \eqref{equ:def of J collective}.
\end{prpstn}
\begin{proof}
  The projection matrix $\Pi_{\nabla g(\bfX)}$ is block diagonal, with blocks $(\Pi_{\nabla g^{\alpha}(\bfX_{\alpha})})_{\alpha \in S}$.
  Hence, by writing $\Expe{\|\nabla u(\bfX)\|_2^2} = \sum_{\alpha \in S} \Expe{\|\nabla_{\alpha} u(\bfX)\|_2^2}$ we can write
  \[
    \calJ(g) 
    = \sum_{\alpha \in S}
    \Expe{\|\nabla_{\alpha} u(\bfX)\|_2^2 - \|\Pi_{\nabla g^{\alpha}(\bfX_{\alpha})}\nabla_{\alpha} u(\bfX)\|_2^2}
    = \sum_{\alpha \in S}
    \Expe{\|\Piperp_{\nabla g^{\alpha}(\bfX_{\alpha})} \nabla_{\alpha} u(\bfX)\|_2^2}.
  \]
  Finally, we obtain the desired result by noting that
  \[
    \Expe{\|\Piperp_{\nabla g^{\alpha}(\bfX_{\alpha})} \nabla_{\alpha} u(\bfX)\|^2}
    = 
    \Expe{\|\Piperp_{\nabla (g^{\alpha}, id^{\alpha^c})(\bfX)} \nabla u(\bfX)\|^2_2}
    = 
    \calJ((g^{\alpha}, id^{\alpha^c})).
  \]
\end{proof}

As in the two groups setting, a consequence of Proposition \ref{prop:poincare bound decomposition multi} is that minimizing $\calJ$ over $\calG_{\bfm} \equiv \times_{\alpha \in S} \calG^{\alpha}_{m_{\alpha}}$ is equivalent to minimizing $\calJ_{\calX_{\alpha}}$ over $\calG_{m_{\alpha}}^{\alpha}$ for all $\alpha \in S$.
In other words, it is equivalent to solving $N$ instances of the collective setting described in \Cref{sec:collective dimension reduction}, where for each instance one would select $\alpha\in S$, set $Y = (\bfX_{\beta})_{\beta \in S \setminus \{\alpha\}}$, and minimize $g^{\alpha} \mapsto \calJ_{\calX_{\alpha}}(g^{\alpha})$.
As a result one may consider leveraging the surrogates $(\calL_{\calX_{\alpha}, m_{\alpha}})_{\alpha \in S}$ from \Cref{subsec:quadratic surrogate to the truncated loss}.
Note however that one would then use a tensorized sample of the form $((\bfx_{\alpha}^{(i_{\alpha})})_{1\leq i_{\alpha} \leq n_{\alpha}})_{\alpha \in S}$ of size $n_s = \prod_{\alpha \in S} n_{\alpha}$, that is exponential in $N$.
Recall also that from Corollary \ref{coro:approximation error quasi optimal multi}, minimizing $g^{\alpha} \mapsto \calE_{\calX_{\alpha}}(g^{\alpha})$ is only nearly equivalent to minimizing $g \mapsto \calE(g)$.


\section{Toward hierarchical formats}
\label{sec:toward hierarchical formats}

In this section, we discuss a generalization of the notion of $\alpha$-rank, see for example \cite[equation 6.12]{hackbuschTensorSpacesNumerical2019}, which we call the $\alpha$-feature-rank.

\begin{dfntn}[feature-rank]
    For $v : \calX \subset \bbR^d \rightarrow \bbR$ and $\alpha \subset D=\{1, \cdots, d\}$, we define the \emph{$\alpha$-feature-rank} of $v$, denoted $\mathrm{rankf}_{\alpha}(v)$, as the smallest integer $r_{\alpha}$ such that 
    \[
        v(\bfx) = f(g(\bfx_{\alpha}), \bfx_{\alpha^c})
    \]
    for some $g:\calX_{\alpha} \rightarrow \bbR^{r_{\alpha}}$ and $f: \bbR^{r_{\alpha}} \times \calX_{\alpha^c}\rightarrow \bbR$.
\end{dfntn}

We can list a few basic properties of the feature-rank.
Firstly, $\mathrm{rankf}_{D}(v)= 1$.
Secondly, for any $\alpha \subset D$, we can write $v(\bfx) = v(id^{\alpha}(\bfx_{\alpha}), \bfx_{\alpha^c})$, thus $\mathrm{rankf}_{\alpha}(v) \leq |\alpha|$.

Now, some important properties of the $\alpha$-rank of multivariate functions are not satisfied by the $\alpha$-feature-rank.
A first property of the $\alpha$-rank is that $\mathrm{rank}_{\alpha}(v) = \mathrm{rank}_{\alpha^c}(v)$, see for example \cite[Lemma 6.20]{hackbuschTensorSpacesNumerical2019}, while this may not be the case for the feature-rank.
A second property of the $\alpha$-rank, which is important for tree-based tensor networks, is \cite[Proposition 9]{nouyHigherorderPrincipalComponent2019}, which states that for any subspace $V_{\alpha} \subset L^2(\calX_{\alpha}, \mu_{\alpha})$, projection onto $V_{\alpha} \otimes L^2(\calX_{\alpha^c}, \mu_{\alpha^c})$ does not increase the $\alpha$-rank, meaning that for any $v\in L^2(\calX, \mu)$,
\[
    \text{rank}_{\alpha} 
    (\calP_{V_{\alpha} \otimes L^2(\calX_{\alpha^c}, \mu_{\alpha^c})} v)
    \leq 
    \text{rank}_{\alpha}(v).
\]
This property was a core ingredient for obtaining near-optimality results when learning tree-based tensor formats with the leaves-to-root algorithm from \cite{nouyHigherorderPrincipalComponent2019}.
The problem here is that our definition of feature-rank does not satisfy this property anymore, as such a projection can increase $\mathrm{rankf}_{\alpha}$.
This is illustrated in the following example.

\begin{xmpl}
Let $\bfX \sim \mu = \calU([-1,1]^3)$ and consider 
\[
    v:\bfx \mapsto (x_1 + x_2) + (x_1 + x_2)^2 x_3.
\]
Since we can write $v(\bfx) = f(x_1 + x_2, x_3)$ for some function $f$, it holds $\mathrm{rankf}_{\alpha}(v) = 1$ for $\alpha = \{1,2\}$.
Firstly, let us consider the subspace $V_{\alpha} = \spanv{\phi_1, \phi_2}$ of $L^2(\calX_{\alpha}, \mu_{\alpha})$ with orthonormal vectors $\phi_1(\bfx_{\alpha}) = \sqrt{3} x_1$ and $\phi_2(\bfx_{\alpha}) = \sqrt{5} x_2^2$.
We then have that 
\[
    (\calP_{V_{\alpha} \otimes L^2(\calX_{\alpha^c}, \mu_{\alpha^c})} v) : \bfx \mapsto x_1 + x_2^2x_3,
\]
thus $\mathrm{rankf}_{\alpha}(\calP_{V_{\alpha} \otimes L^2(\calX_{\alpha^c}, \mu_{\alpha^c})} v) = 2$.
Let us also consider $W_{\alpha} = \Sigma(x_1\mapsto x_1) \otimes \Sigma(x_2 \mapsto x_2^2)$.
We then have that 
\[
    (\calP_{W_{\alpha} \otimes L^2(\calX_{\alpha^c}, \mu_{\alpha^c})} v) : \bfx \mapsto x_1 + (x_1^2 + x_2^2) x_3,
\]
thus $\mathrm{rankf}_{\alpha}(\calP_{W_{\alpha} \otimes L^2(\calX_{\alpha^c}, \mu_{\alpha^c})} v) = 2$.
As a result, for both examples $V_{\alpha}$ and $W_{\alpha}$, projection increased the $\alpha$-feature-rank.
\end{xmpl}


\section{Numerical experiments}
\label{sec:numerical experiments structured dim red}

\subsection{Setting}

In this section we apply the collective dimension reduction approach described in \Cref{sec:collective dimension reduction} to a polynomial of $\bfX \sim \calU(\calX)$ with $\calX = (-1,1)^d$ and $d=8$, with coefficients depending on $Y \sim \calU(\calY)$ with $\calY = (-1,1)$, where $\bfX$ and $Y$ are independent.
For $a \geq 1$ we define $u_a$ by
\begin{equation}
    u_a(\bfx, y) := 
    \sum_{k=1}^{a} (\bfx^T Q_k \bfx)^2
    \sin(\frac{\pi k}{2 a} y), 
\end{equation}
with symmetric matrices $Q_k := \frac{1}{2}(\mathbbm{1}_{i-j = k-1} + \mathbbm{1}_{j-i = k-1})_{ij} \in \bbR^{d\times d}$ for $1\leq k\leq a$.
In this context, we can express $u(\bfX, Y)$ as a function of degree-$2$ polynomial features in $\bfX$, as we can write $u_a(\bfX, Y) = f(g(\bfX), Y)$ with $g(\bfx)= (\bfx^T Q_k \bfx)_{1\leq k\leq a}$ and with $f(\bfz, y) = \sum_{1\leq k\leq a} z_k^2 \sin(\frac{\pi k}{2 a} y)$.
We consider two cases.
Firstly $a=m=3$, secondly $a=3$ and $m=2$.
In the first case $u_a$ can be exactly represented as a function of $m$ degree $2$ polynomial features, but not in the other case.

In our experiments, we will monitor $4$ quantities.
The first two are the Poincar\'e inequality based quantity $\calJ_{\calX}(g)$ defined in \eqref{equ:def of J collective} and the final approximation error $e_g(f)$ defined by 
\[
    e_g(f) := \Expe{|u(\bfX, Y) - f(g(\bfX), Y) |^2}^{1/2}.
\]
We estimate these quantities with their Monte-Carlo estimators on a test sample $\Xi^{test} \subset \calX \times \calY$ of size $N^{test}=1000$, not used for learning.
We also monitor the Monte-Carlo estimators $\widehat{\calJ}_{\calX}(g)$ and $\widehat{e}_g(f)$ on some training set $\Xi^{train}\subset\calX \times  \calY$ of various sizes $N^{train}$, which will be the quantities directly minimized to compute $g$ and $f$.
More precisely, we draw $20$ realizations of $\Xi^{train}$ and $\Xi^{test}$ and monitor the quantiles of those $4$ quantities over those $20$ realizations.

We consider feature maps of the form \eqref{equ:Gm vector space sphere} with $\Phi : \bbR^d \rightarrow \bbR^K$ a multivariate polynomial of total degree at most $\ell+1 = 2$, excluding the constant polynomial so that $\mathrm{rank}(\nabla \Phi(\bfX)) =d$ almost surely.
Note also that this definition ensures that $id \in \spanv{\Phi_1, \cdots, \Phi_K}$, thus $\calG_m$ contains all linear feature maps, including the one corresponding to the active subspace method.

We compare two procedures for constructing the feature map.
The first procedure, which we consider as the reference, is based on a preconditioned nonlinear conjugate gradient algorithm on the Grassmann manifold $\mathrm{Grass}(m, K)$ to minimize $G \mapsto \widehat \calJ_{\calX}(G^T \Phi)$.
For this procedure, the training set 
\[
    \Xi^{train} = (\bfx^{(k)}, y^{(k)})_{1\leq k\leq N^{train}}
\]
is drawn as $N^{train}$ samples of $(\bfX, Y)$ using a Latin hypercube sampling method.
We use $\hat{\Sigma}(G) \in \bbR^{K\times m}$ as a preconditioning matrix at the point $G\in\bbR^{K\times m}$, which is the Monte-Carlo estimation of $\Sigma(G)$ defined in \cite[Proposition 3.2]{bigoniNonlinearDimensionReduction2022}.
We choose as initial point the matrix $G^0\in\bbR^{K\times m}$ which minimizes $\widehat{\calJ}$ on the set of linear features, which corresponds to the active subspace method.
We denote this reference procedure as GLI, standing for Grassmann Linear Initialization.

The second procedure consists of taking the feature map that solves Proposition \ref{prop:L collective as generalized ev}, with $H_{\calX, m}$ replaced with its Monte-Carlo estimator on the tensorized set 
\[
    \Xi^{train} = (\bfx^{(i)}, y^{(j)})_{1\leq i\leq n_{\calX},1\leq j\leq n_{\calY}}
\] 
of size $N^{train} = n_{\calX} n_{\calY}$ with $n_{\calY}=5$ fixed.
The two samples $(\bfx^{(i)})_{1\leq i\leq n_{\calX}}$ and $(y^{(j)})_{1\leq j\leq  n_{\calY}}$ are independent samples drawn using Latin hypercube sampling.
Estimating $H_{\calX, m}$ includes estimating $M(\bfx^{(i)})$ and $V_m(\bfx^{(i)})$ with their Monte-Carlo estimators on $(y^{(j)})_{1\leq j\leq n_{\calY}}$ for all $1\leq i\leq n_{\calX}$.
Note that $R = \Expe{\nabla \Phi(\bfX)^T \nabla \Phi(\bfX)}$ is exactly computed thanks to the choice of $\Phi$.
We denote this procedure as SUR, standing for SURrogate.
We emphasize the fact that the methods SUR and GLI are not performed on the same training sets, although the sizes of the training sets are the same.

Once $g\in\calG_m$ is learnt, we then perform a classical regression task to learn a regression function $f:\bbR^m \times \bbR \rightarrow \bbR$, with $g(\bfX)\in\bbR^m$ and $Y \in \bbR$ as input variables and $u(\bfX, Y)\in\bbR$ as the output variable.
In particular here, we have chosen to use kernel ridge regression with a Gaussian kernel $\kappa(\bft^{(1)}, \bft^{(2)}) := \exp(-\gamma \|\bft^{(1)} - \bft^{(2)}\|^2_{2})$ for any $\bft^{(1)},\bft^{(2)}\in\bbR^m$ and some hyperparameter $\gamma>0$.
Then with $\{\bfz^{(i)}\}_{1\leq i\leq N^{train}} := \{(g(\bfx), y) : (\bfx, y) \in \Xi^{train}\}$, we consider
\[
    f: \bfz  \mapsto \sum_{i=1}^{N^{train}} a_i \kappa(\bfz^{(i)}, \bfz),      
\]
with $\bfa :=  ( B + \alpha I_{N^{train}})^{-1} \bfu \in \bbR^{N^{train}}$ for some regularization parameter $\alpha >0$, where $B := (\kappa(\bfz^{(i)}, \bfz^{(j)}))_{1\leq i,j\leq N^{train}}$ and $\bfu := u(\Xi^{train}) \in \bbR^{N^{train}}$.
Here the kernel parameter $\gamma$ and the regularization parameter $\alpha$ are hyperparameters learnt using a $10$-fold cross-validation procedure, such that $\log_{10}(\gamma)$ is selected from $30$ points uniformly spaced in $[-6, -2]$, and $\log_{10}(\alpha)$ is selected from $40$ points uniformly spaced in $[-11, -5]$.
Note that these sets of hyperparameters have been chosen arbitrarily to ensure a compromise between computational cost and flexibility of the regression model.
Note also that with additional regularity assumptions on the conditional expectation $(\bfz, y) \mapsto \Expe{u(\bfX, Y)| (g(\bfX), Y)=(\bfz, y)}$ it may be interesting to consider a Mat\'ern kernel instead of the Gaussian kernel.

The cross-validation procedures as well as the Kernel ridge regression rely on the library \emph{sklearn} \cite{pedregosaScikitlearnMachineLearning2011}.
The optimization on Grassmann manifolds relies on the library \emph{pymanopt} \cite{JMLR:v17:16-177}.
The orthonormal polynomial feature maps rely on the Python library \emph{tensap} \cite{nouyAnthonynouyTensap2020}.
The code underlying this work is freely available as a Zenodo deposit \cite{pascoAlexandrepascoTensapV171poincarelearningcollectivepaper2026}.

\subsection{Results and observations}

\begin{figure}[h]
    \centering
    \includegraphics[page=2, width=0.99\textwidth]{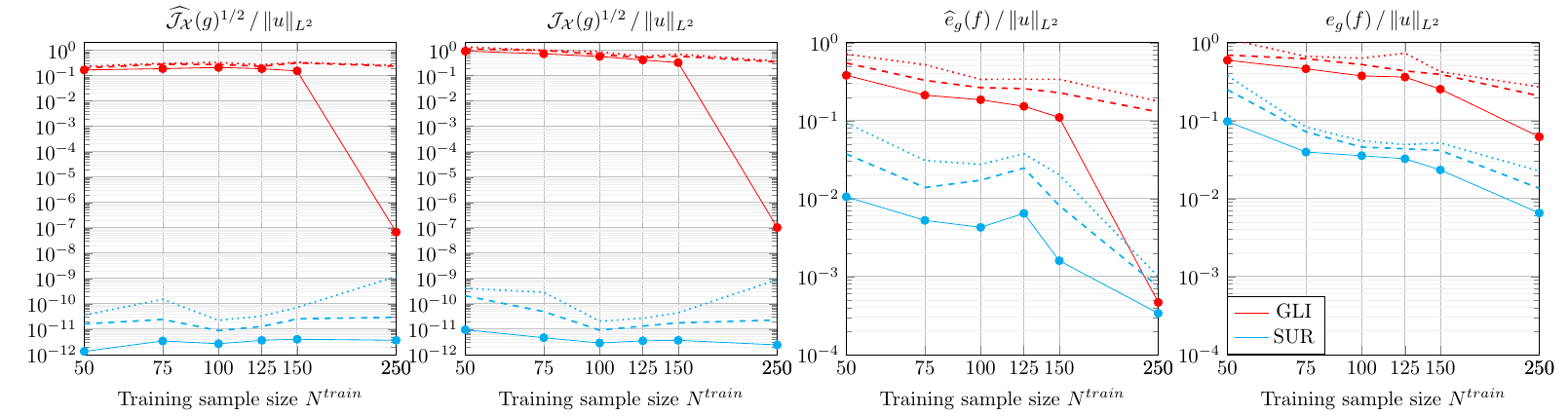}
    \caption[caption]{\footnotesize 
    Evolution of quantiles with respect to the size of the training sample for $u_3$ with $m=3$.
    The quantiles $50\%$, $90\%$ and $100\%$ are represented respectively by the continuous, dashed and dotted lines.
    }
    \label{fig:u3 m3 vs ntrain}
\end{figure}

Let us start with $u_3$ approximated with $m=3$ features and $a=3$, for which results are displayed in \Cref{fig:u3 m3 vs ntrain}.
Firstly, for all values of $N^{train}$, we observe that SUR always yields the minimizer of $\widehat \calJ_{\calX}$, which turns out to be $0$ as for $\calJ_{\calX}$.
On the other hand, GLI mostly fails to achieve such a result for $N^{train}\leq 150$, and sometimes fails to achieve such a result for $N^{train} = 250$.
A large performance gap is also observed regarding $\widehat{e}_g(f)$ and $e_g(f)$.
We also observe that, although the minimum of $e_g$ over all measurable functions should be $0$, its minimum over the chosen regression class is not $0$.

\begin{figure}[h]
    \centering
    \includegraphics[page=1, width=0.99\textwidth]{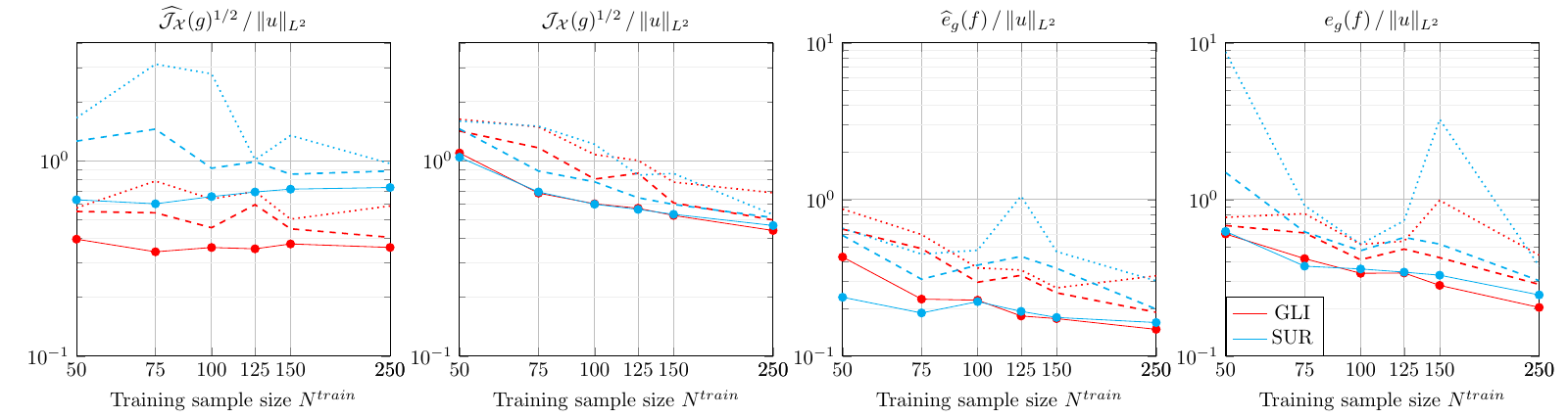}
    \caption[caption]{\footnotesize 
    Evolution of quantiles with respect to the size of the training sample for $u_3$ with $m=2$.
    The quantiles $50\%$, $90\%$ and $100\%$ are represented respectively by the continuous, dashed and dotted lines.
    }
    \label{fig:u3 m2 vs ntrain}
\end{figure}

Let us continue with $u_3$ approximated with $m=2 < a$ features, for which results are displayed in \Cref{fig:u3 m2 vs ntrain}.
We first observe that GLI performs better at minimizing $\widehat{\calJ}_{\calX}$ than SUR, although the corresponding performance on $\calJ_{\calX}$ is rather similar.
We then observe that SUR and GLI perform mostly similarly regarding the regression errors $\widehat{e}_g(f)$ and $e_g(f)$.
However, SUR suffers from important performance gaps between $\widehat{e}_g(f)$ and $e_g(f)$ in some worst-case errors.
This might be due to the small size $n_{\calY}=5$ for the sample of $Y$.

\section{Conclusion and perspectives}
\label{sec:conclusion structured features}

\subsection{Conclusion}

In this paper we analyzed two types of nonlinear dimension reduction problems in a regression framework.

We first considered a collective dimension setting, which consists in learning a feature map suitable to a family of functions.
Considering Poincar\'e inequality based methods, we extended the surrogate approach developed in \cite{nouySurrogatePoincareInequalities2026} to the collective setting.
We showed that for polynomial feature maps, and under some assumptions, our surrogate can be used as an upper bound of the Poincar\'e inequality based loss function.
Moreover, the surrogate we introduced is quadratic with respect to the feature maps, thus well suited for optimization procedures.
In particular when the features are taken from a finite dimensional linear space, then minimizing the surrogate is equivalent to finding the eigenvectors associated to the smallest generalized eigenvalues of some matrix pencil.
The main practical limitation of our surrogate is that it cannot be used with arbitrary samples, as it requires a specific sampling structure such as tensorized samples.

We then considered a two groups setting, which consists in learning two different feature maps associated to disjoint groups of input variables.
We drew the parallel with functional singular value decomposition, pointing out the main similarities and differences.
We also considered a multiple groups setting, which consists in separating the input variables into more than two groups and learning corresponding feature maps.
We drew the parallel with the Tucker tensor format, which allowed us to obtain a near-optimality result similar to the near-optimality of the higher order singular value decomposition.
More precisely, the multiple groups setting is almost equivalent to several instances of the collective setting.
Additionally, when considering Poincar\'e inequality based methods, the equivalence holds.
We also discussed extending the analysis towards hierarchical format, trying to draw the parallel with tree-based tensor networks.
However, we were only able to draw some pessimistic results.
In particular, we investigated a new notion of rank, which unfortunately  lacks some important properties leveraged in the analysis of tree-based tensor networks.

Finally, we illustrated our surrogate method in the collective setting on a numerical example.
We observed that when a representation with low dimensional features existed, our method successfully identified them, while direct methods for minimizing the Poincar\'e inequality based loss function mostly failed.
However, we observed that when such a low-dimensional representation did not exist, our method performed mostly similarly to the other one.
In particular in the worst-case scenario we observed that the regression procedure in our approach may be very challenging, which is probably due to the tensorized sampling strategy.

\subsection{Perspectives}

Let us mention three main avenues for future work.
The first perspective is to find intermediate regimes of interest for the class of regression functions.
Indeed, in this paper we discussed only the linear and measurable cases, which essentially constitute two extreme opposites of the possible choices of classes of feature maps.
The second perspective is to further investigate the fundamental properties of the collective, two-group and multiple-group settings.
Indeed, in our analysis we showed near-optimality results assuming that these problems do admit solutions, which we have not properly demonstrated.
The third perspective is to extend our surrogate approach to the Bayesian inverse problem setting.
Indeed, recent works leveraged gradient-based functional inequalities to derive certified nonlinear dimension reduction methods for approximating the posterior distribution in this framework \cite{liPrincipalFeatureDetection2024,liSharpDetectionLowdimensional2025}.
Extending our surrogate methods may improve the learning procedure of nonlinear features in such a setting.


\section*{Tool and computational resource disclosure}
\label{sec:acknowledgments}

In addition to the Python libraries mentioned in \Cref{sec:numerical experiments structured dim red}, Artificial Intelligence (AI) assistance was used to generate Figures \ref{fig:illustration dim red non collective} and \ref{fig:illustration dim red collective}, and to generate the notations list in Appendix \ref{apx:notations}.

\section*{Data Availability Statement}
\label{sec:data availability}

The code used in this paper is available online as a Zenodo deposit \cite{pascoAlexandrepascoTensapV171poincarelearningcollectivepaper2026}.


\bibliographystyle{amsplain}  
\bibliography{main}

\appendix

\section{Poincar\'e inequalities for dimension reduction}
\label{apx:poincare inequalities for dimension reduction}

\subsection{Poincar\'e inequalities for measures}

A well known result is the Gaussian Poincar\'e inequality, see for example \cite[Theorem 3.20]{boucheronConcentrationInequalitiesNonasymptotic2013}.
It states that if $\bfX\in\bbR^d$ is a vector whose entries are i.i.d. centered standard Gaussian random variables, then for any real valued continuously differentiable function $h : \bbR^d \rightarrow \bbR$ it holds
\[
    \Expe{\left(h(\bfX) - \Expe{h(\bfX)} \right)^2}
    \leq
    \Expe{\|\nabla h(\bfX)\|_2^2},
\]
where the left-hand side is the variance of $h(\bfX)$.
More generally, we say that the probability measure $\mu_{\calX}$ of a random variable $\bfX$ satisfies a Poincar\'e inequality if there exists a constant $C(\mu_{\calX}) < +\infty$ such that for every continuously differentiable function $h : \calX \rightarrow \bbR$ it holds
\begin{equation}
\label{equ:poincare inequality general measure}
    \Expe{\left(h(\bfX) - \Expe{h(\bfX)} \right)^2}
    \leq
    C(\mu_{\calX})\Expe{\|\nabla h(\bfX)\|_2^2}.
\end{equation}
The smallest constant for which the above inequality holds is called the Poincar\'e constant associated to the measure $\mu_{\calX}$.
Note that such an inequality can only hold when the support $\calX$ of the measure is connected, otherwise $C(\mu_{\calX}) = +\infty$.
Typical examples with explicit bounds on the Poincaré constant are the uniform measures on a bounded convex set, for which $C(\mu_{\calX}) \leq \mathrm{diam}(\calX)^2 / \pi^2$ and the non-standard Gaussian measure with covariance $\Sigma$, for which $C(\mu_{\calX}) = \lambda_{\max}(\Sigma)$.
Both measures are examples of $s$-concave measures, as defined in \ref{def:s-concave measures}, the first with $s \leq 1/\dim(\calX)$ and the second with $s \leq 0$.
More generally, \cite[Theorem 1.2]{bobkovWeightedIsoperimetricPoincaretype2009} implies that for any $s\in[0, 1/\dim(\calX)]$, an $s$-concave probability measure satisfies \eqref{equ:poincare inequality general measure} with constant $C(\mu_{\calX}) \leq c_0 \exp\Expe{\log\|\bfX\|_2^2}$, where $c_0<\infty$ is a universal constant.
Note that for $s<0$, which includes probability measures with heavy tails such as the Cauchy distribution, only a \emph{weighted} Poincar\'e inequality is guaranteed to hold, which is weaker.

The probability measures mentioned previously are all absolutely continuous with respect to the Lebesgue measure, but \eqref{equ:poincare inequality general measure} can also hold for measures supported on a more general smooth submanifold $\calX$ of $\bbR^d$, endowed with the Euclidean metric. 
The gradient on the right-hand side of \eqref{equ:poincare inequality general measure} is then read as the \emph{Riemannian} gradient.
We refer the reader to \cite{absilOptimizationAlgorithmsMatrix2008} for an introduction to Riemannian geometry.
For example, if $\mu_{\calX}$ is the uniform measure on the sphere $\calX=\mathbb{S}^{d-1}$, then $C(\mu_{\calX}) = 1/(d-1)$.
More generally, if $\mu_{\calX}$ is the uniform measure on a compact connected smooth Riemannian manifold $\calX$, then \cite[Theorem 2.10]{hebeyNonlinearAnalysisManifolds2000} states that $C(\mu_{\calX}) < \infty$.
Note that the Poincar\'e constant can also be defined as the inverse of the smallest nonzero eigenvalue of the (possibly weighted) Laplacian operator generalized to nonlinear manifolds.

\subsection{Stronger inequality required for dimension reduction}

In order to obtain the inequalities \eqref{equ:upper bound approx error linear} or \eqref{equ:upper bound approx error nonlinear}, the Poincar\'e inequality satisfied by $\mu_{\calX}$ is not sufficient.
Instead, we need that for all feature maps $g \in \calG_m$ and all features $\bfz \in g(\calX)$, the conditional measure $\mu_{\calX | g(\bfX)=\bfz}$ satisfies a Poincar\'e inequality with a finite constant independent of $g$ and $\bfz$,
\begin{equation}
\label{equ:def uniform poincare constant}
    C(\mu_{\calX},\calG_m) := \sup_{g\in\calG_m} \sup_{\bfz \in g(\calX)} C(\mu_{\calX | g(\bfX)=\bfz}).
\end{equation}
Note that the constant $C$ in \eqref{equ:upper bound approx error nonlinear} can be chosen as $C(\mu_{\calX},\calG_m)$.
More precisely, one can obtain \eqref{equ:upper bound approx error nonlinear} by first conditioning on $g(\bfX)=\bfz$, then applying \eqref{equ:poincare inequality general measure} to the conditional measure $\mu_{\calX | g(\bfX)=\bfz}$ and the function $h_{\bfz} = u_{|g^{-1}(\{\bfz\})} : g^{-1}(\{\bfz\}) \rightarrow \bbR$, which denotes the restriction of $u$ to the level set $g^{-1}(\{\bfz\})$, then by taking the expectation over $\bfz$.
The last step consists in expressing the Riemannian gradient of $h_{\bfz}$, which simplifies to \eqref{equ:upper bound approx error linear} if $g$ is linear.

For linear feature maps, the level sets of $g$ are affine subspaces of $\bbR^d$.
This is useful because some properties of measures that are sufficient conditions for a Poincar\'e inequality to hold are preserved under conditioning on an affine space.
A very favorable case is when $\mu_{\calX}$ is the standard Gaussian distribution.
Then, the conditional measures are also Gaussian distributions, whose maximum directional variance is bounded by $1$, which yields $C(\mu_{\calX},\calG_m) = 1$ as stated in \cite{zahmGradientBasedDimensionReduction2020}.
Other examples of such properties are uniformity on a convex set and strong log-concavity, and explicit upper bounds on $C(\mu_{\calX},\calG_m)$ are known, see \cite{parenteGeneralizedBoundsActive2020,zahmCertifiedDimensionReduction2022}.
We briefly discuss possible extensions to $s$-concavity in Remark \ref{rmk:extension to s concave}.

For nonlinear feature maps, the level sets of $g$ may not be affine subspaces anymore, but rather smooth submanifolds of $\bbR^d$, under a full rank assumption on the Jacobian of $g$.
Controlling the Poincar\'e constant associated to all the conditional measures is then a challenging task, as their support may have various forms.
The approach considered in \cite{bigoniNonlinearDimensionReduction2022,romorKernelbasedActiveSubspaces2022,nouySurrogatePoincareInequalities2026} consists in essentially ignoring $C(\mu_{\calX}, \calG_m)$ while focusing on the minimization of $\calJ$ on a set of classical feature maps such as polynomials or trigonometric functions.
In fact, these choices do not guarantee the connectivity of the level sets of $g$, which may lead to $C(\mu_{\calX},\calG_m) = +\infty$.
Still, empirical evidence shows that this can be a relevant approach, although performance in large-scale applications has not been demonstrated yet.
On the other hand, the approach considered in \cite{zhangLearningNonlinearLevel2019,verdiereDiffeomorphismbasedFeatureLearning2025} uses feature maps based on diffeomorphisms, which helps ensure the connectivity of the level sets, although it does not bound $C(\mu_{\calX}, \calG_m)$.
In particular, the authors in \cite{verdiereDiffeomorphismbasedFeatureLearning2025} proposed a heuristic penalization of $\calJ(g)$ by some Kullback-Leibler divergence to the standard Gaussian measure.
Let us also mention \cite{pillaud-vivienStatisticalEstimationPoincare2020} in which, motivated by \cite{lelievreGeneralTwoscaleCriteria2009}, the authors proposed to maximize $C(g_{\#\mu_{\calX}})$ as a heuristic for minimizing $C(\mu_{\calX | g(\bfX)=\bfz})$.
They then leveraged the fact that the former can be estimated from a given sample, contrary to the latter.
Although the setting was different, it may be relevant to consider such a heuristic in our dimension reduction setting.

\begin{rmrk}
\label{rmk:extension to s concave}
    Using Definition \ref{def:s-concave measures} one should be able to show that $s$-concavity is preserved when conditioning on linear features.
    If so, then from \cite{bobkovWeightedIsoperimetricPoincaretype2009}, one would still need to bound the conditional geometric mean $\exp \Expe{\log\|\bfX\|_2 | g(\bfX)=\bfz}$, up to some (conditional) shift as the Poincar\'e inequality is shift invariant.
    This is left to further investigation.
\end{rmrk}

\subsection{On the minimization of the loss function}

In this section we briefly discuss the minimization of the loss function
\[
    \calJ(g) = \Expe{\|\Piperp_{\nabla g(\bfX)} \nabla u(\bfX)\|_2^2}
    = \Expe{\|\nabla u(\bfX)\|_2^2} - \Expe{\|\Pi_{\nabla g(\bfX)} \nabla u(\bfX)\|_2^2},
\]
which we defined in \eqref{equ:upper bound approx error nonlinear} and was originally introduced in \cite{bigoniNonlinearDimensionReduction2022}.
We refer to \cite{bigoniNonlinearDimensionReduction2022,verdiereDiffeomorphismbasedFeatureLearning2025,nouySurrogatePoincareInequalities2026} for a broader discussion.
The above expression can be interpreted as the mean projection error of $\nabla u(\bfX)$ onto the subspace $\spanv{\nabla g(\bfX)}$.

\begin{figure}[h]
    \centering
    \includegraphics[page=1, width=0.49\textwidth]{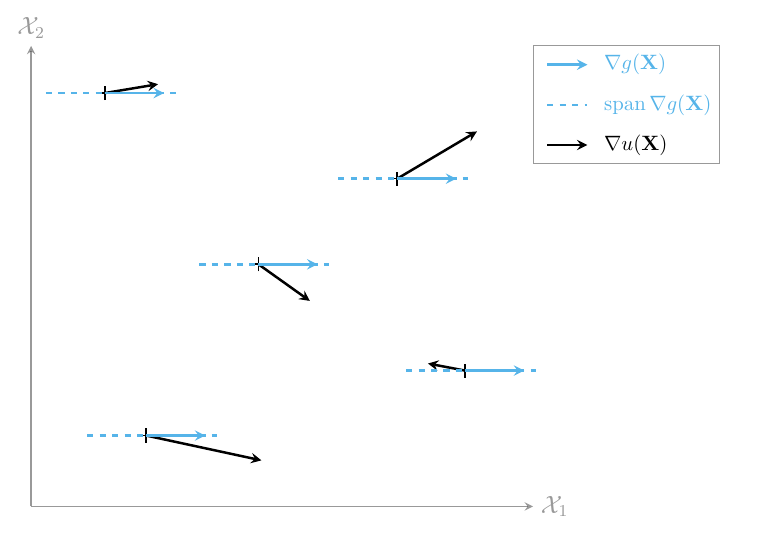}
    \includegraphics[page=2, width=0.49\textwidth]{./figures/drawing/drawing.pdf}
    \caption[]{\footnotesize 
        Illustration of gradient alignment in a non-collective setting with $m=1$ linear features (left) and nonlinear features (right), for five sample points of $\bfX\in\bbR^2$.
    }
    \label{fig:illustration dim red non collective}
\end{figure}

If $g$ is linear, of the form $g(\bfx) = G^T\bfx$ for some $G\in \bbR^{d\times m}$, then $\nabla g(\bfX) = G$.
Importantly, $\spanv{\nabla g(\bfX)}$ is constant with respect to $\bfX$, as illustrated in the left part of Figure \ref{fig:illustration dim red non collective}.
Then, $\calJ(g)$ is the mean squared projection error of $\nabla u(\bfX)$ onto the subspace $\spanv{G}$.
The optimal subspace in this sense is commonly obtained via (uncentered) principal component analysis (PCA) of the random vector $\nabla u(\bfX)$, and equals the dominant eigenspace of the second moment matrix $\Expe{\nabla u(\bfX) \nabla u(\bfX)^T}$.
This setting is illustrated in the left part of Figure \ref{fig:illustration dim red non collective}.
If $g$ is nonlinear, then $\spanv{\nabla g(\bfX)}$ may no longer be constant.
This constitutes one of the reasons why minimizing $\calJ$ on a set of nonlinear features may be challenging.
This setting is illustrated in the right part of Figure \ref{fig:illustration dim red non collective}.

\subsection{Collective setting}

As explained above, we mentioned that \eqref{equ:upper bound approx error nonlinear} was obtained by applying \eqref{equ:poincare inequality general measure} to the conditional measure $\mu_{\calX | g(\bfX)=\bfz}$ and the function $u_{|g^{-1}(\{\bfz\})}$.
The collective setting is very similar, as \eqref{equ:E bounded by J collective} is obtained by applying \eqref{equ:poincare inequality general measure} to the functions $(u_y)_{|g^{-1}(\{\bfz\})} : g^{-1}(\{\bfz\}) \rightarrow \bbR$ for all $y\in\calY$.
One difference between the collective and the non-collective settings is that a single function can always be represented as a function of a single feature, for example by taking the function itself as the feature map, while a family of functions may intrinsically require more features.
For example, it requires at least two features in $\bfX$ to exactly represent $u_Y(\bfX) = \sin(X_1) + \cos(Y X_2)$.

\begin{figure}[h]
    \centering
    \includegraphics[page=3, width=0.49\textwidth]{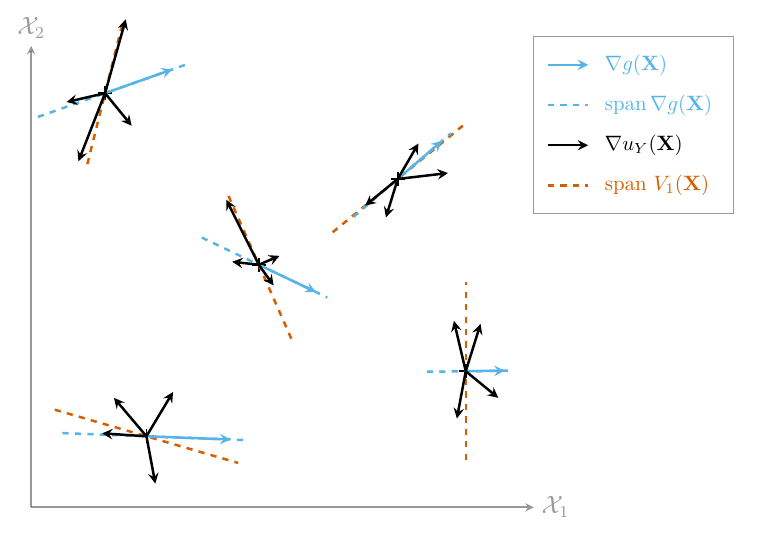}
    \caption[]{\footnotesize 
        Illustration of gradient alignment in the collective setting with $m=1$ nonlinear feature for five sample points of $\bfX\in\bbR^2$ and four sample points of $Y$.
    }
    \label{fig:illustration dim red collective}
\end{figure}

This can be seen in the case $m=1$ when comparing Figures \ref{fig:illustration dim red non collective} and \ref{fig:illustration dim red collective}.
In the latter, we see that a single direction $\spanv{\nabla g(\bfx)}$ may not be able to perfectly align with all the vectors $\nabla u_Y (\bfx)$.
More rigorously, this is a consequence of \eqref{equ:J collective lower bound eigen values}, in which we showed that 
\[
    \calJ_{\calX}(g) 
    = \Expe{\|\Piperp_{\nabla g(\bfX)} \nabla u_Y(\bfX)\|_2^2}
\]
is lower bounded by $\varepsilon_m$ defined in \eqref{equ:J collective lower bound eigen values}, which is independent of the choice of the feature maps.
To understand why, recall that $\varepsilon_m$ is defined as the mean squared projection error of $\nabla u_Y(\bfX)$ onto $\spanv{V_m(\bfX)}$, the latter corresponding to the dotted orange line in Figure \ref{fig:illustration dim red collective}.
Recall also that $V_m(\bfx)$ is defined via (uncentered) PCA of the random vector $\nabla u_Y(\bfx)$, where $\bfx$ is fixed, and thus spans the subspace that minimizes the mean squared projection error conditional on $\bfX=\bfx$.
As a result, the best we can hope for is that $\spanv{V_m(\bfx)}$ is well aligned with $\spanv{\nabla g(\bfx)}$, the dotted blue line in Figure \ref{fig:illustration dim red collective}.
However, the conditions on $u$ under which there exists $g:\bbR^d \rightarrow \bbR^m$ such that $\spanv{\nabla g(\bfX)} = \spanv{V_m(\bfX)}$, almost surely, remain unclear to us.

The underlying idea behind our approach from \Cref{subsec:quadratic surrogate to the truncated loss} is that instead of minimizing the projection error of $V_m(\bfX)$ onto $\spanv{\nabla g(\bfX)}$, we minimize a (weighted) projection error of $\nabla g(\bfX)$ onto $\spanv{V_m(\bfX)}$,
\[
    \calL_{\calX,m}(g) =
    \Expe{
        \lambda_1(M(\bfX))
        \|\Piperp_{V_m(\bfX)} \nabla g(\bfX)\|_F^2
    },
\]
with suitable normalization constraints on $g$.
Note that Figure \ref{fig:illustration dim red collective} also illustrates why our approach requires multiple realizations of $Y$ at each sample point of $\bfX$, which is the main limitation of our approach.

\section{Notations}
\label{apx:notations}

\newlength{\notationlabelwidth}
\setlength{\notationlabelwidth}{0.20\textwidth}
\newcommand{\notationlistlabel}[1]{#1\hfil}
\newenvironment{notationlist}
{%
    \begin{list}{}{%
        \setlength{\labelwidth}{\notationlabelwidth}%
        \setlength{\labelsep}{1.25em}%
        \setlength{\leftmargin}{\dimexpr\labelwidth+\labelsep\relax}%
        \setlength{\itemindent}{0pt}%
        \setlength{\listparindent}{0pt}%
        \let\makelabel\notationlistlabel
    }%
    \raggedright
}
{%
    \end{list}
}

This appendix gathers the main notation used throughout the manuscript.  
Unless specified otherwise, vectors are column vectors.  
The references below point to the first definition, or to the first passage fixing the notation when the definition is not displayed.

\subsection*{General dimension-reduction setting}

\begin{notationlist}
    \item[$d$, $m$]
    The input dimension and the number of learned features, respectively, with $m\leq d$; see \Cref{sec:structured features introduction}.

    \item[$\calX$, $\calY$]
    The state spaces of the input variables $\bfX\in\calX\subset\bbR^d$ and $Y\in\calY$; see \Cref{sec:structured features introduction}.

    \item[$\mu_{\calX}$, $\mu_{\calY}$]
    The probability distributions of $\bfX$ and $Y$.  
    In the collective setting they are independent, so the joint distribution is $\mu_{\calX}\otimes\mu_{\calY}$; see \Cref{sec:structured features introduction}.

    \item[$u$, $u_y$]
    The target function $u:\calX\times\calY\to\bbR$ in the collective setting, with $u_y:=u(\cdot,y)$; see the opening of \Cref{sec:collective dimension reduction}.  
    In the grouped setting, $u:\calX\to\bbR$; see \Cref{sec:two variables approach}.

    \item[$g$, $f$, $\hat u$]
    In the non-collective setting, $g:\calX\to\bbR^m$, $f:\bbR^m\to\bbR$, and $\hat u(\bfx)=f(g(\bfx))$. 
    In the collective setting, the same feature map $g:\calX\to\bbR^m$ is shared by every function $u_y$, while $f:\bbR^m\times\calY\to\bbR$ and $\hat u(\bfx,y)=f(g(\bfx),y)$; see \Cref{sec:structured features introduction}.

    \item[$\calG_m$, $\calF_m$]
    The admissible classes of $m$-dimensional feature maps and profile functions, respectively; see the opening of \Cref{sec:collective dimension reduction}.

    \item[$\nabla g(\bfx)$]
    The transposed Jacobian of $g$, whose $j$th column is
    $\nabla g_j(\bfx)$, so that $\nabla g(\bfx)\in\bbR^{d\times m}$; see the paragraph following \eqref{equ:upper bound approx error nonlinear}.

    \item[$\spanv{A}$, $\Pi_A$, $\Piperp_A$]
    The column span of a matrix $A$, the orthogonal projector onto that span, and its complementary projector $\Piperp_A:=I-\Pi_A$, respectively.
    The projector notation is introduced in \eqref{equ:upper bound approx error linear}, and the complementary projector in \eqref{equ:def princial components gradient u}.

    \item[{$\Proba[Z]{\cdot}$, $\Expe[Z]{\cdot}$}]
    Probability and Expectation with respect to $Z$ (the subscript is omitted for the joint expectation), respectively.

    \item[$\lambda_i(A)$, $\sigma_i(A)$]
    The $i$th eigenvalue of a symmetric matrix $A$, ordered decreasingly, and the $i$th singular value of a matrix $A$, also ordered decreasingly.
    These conventions are fixed after \eqref{equ:def princial components gradient u} and in Lemma \ref{lem:bounds on J collective truncated}, respectively.

    \item[$\calP_V$, $\calP_V^\perp$]
    The orthogonal projector onto a closed subspace $V$ of an $L^2$ space and its complementary projector; see \eqref{equ:def svd} and \eqref{equ:approximation error decomposition multi}.
\end{notationlist}

\subsection*{Collective losses and their surrogate}

\begin{notationlist}
    \item[$\calE_{\calX}(g)$]
    The smallest mean squared reconstruction error obtainable from $g(\bfX)$ and $Y$ when the profile function is only required to be measurable; see \eqref{equ:def reconstruction error collective}.

    \item[$\calJ(g)$, $\calJ_{\calX}(g)$]
    The Poincar\'e-inequality-based loss functions in the non-collective and collective settings, respectively; see \eqref{equ:upper bound approx error nonlinear} and \eqref{equ:def of J collective}.

    \item[$C(\mu)$, $C(\mu,\calG_m)$]
    The Poincar\'e constant of a measure $\mu$ and the uniform Poincar\'e constant over conditional measures induced by maps in $\calG_m$; see \eqref{equ:poincare inequality general measure} and \eqref{equ:def uniform poincare constant}.

    \item[$M(\bfx)$, $V_m(\bfx)$]
    The conditional uncentered second-moment matrix of $\nabla u_Y(\bfx)$ and a matrix with orthonormal columns spanning its dominant $m$-dimensional eigenspace, respectively; see \eqref{equ:def covariance matrix gradient u} and \eqref{equ:def princial components gradient u}.

    \item[$\varepsilon_m$]
    A lower bound on $\calJ_{\calX}$, defined as the mean-squared projection error of $\nabla u_Y(\bfX)$ onto $\spanv{V_m(\bfX)}$; see \eqref{equ:J collective lower bound eigen values}.

    \item[$\calJ_{\calX,m}(g)$]
    The truncated collective loss function, obtained by replacing $\nabla u_Y(\bfX)$ with $\Pi_{V_m(\bfX)}\nabla u_Y(\bfX)$ in the definition of $\calJ_{\calX}(g)$; see \eqref{equ:def J collective truncated}.

    \item[$\calL_{\calX,m}(g)$]
    The quadratic loss function used as a surrogate to $\calJ_{\calX, m}$; see \eqref{equ:def L collective}.

    \item[$\Phi$, $G$, $K$]
    A fixed dictionary $\Phi:\calX\to\bbR^K$, its coefficient matrix
    $G\in\bbR^{K\times m}$, and the dictionary dimension, used to parametrize
    $g(\bfx)=G^T\Phi(\bfx)$; see \eqref{equ:Gm vector space sphere}.

    \item[$R$]
    The feature-gradient Gram matrix $R=\Expe{\nabla\Phi(\bfX)^T\nabla\Phi(\bfX)}$, used to impose $G^TRG=I_m$; see \eqref{equ:Gm vector space sphere} and the paragraph preceding Proposition~\ref{prop:L collective as generalized ev}.

    \item[$H_{\calX,m}$]
    The positive semidefinite matrix representing the surrogate through $\calL_{\calX,m}(G^T\Phi)=\Trace{G^TH_{\calX,m}G}$; see \eqref{equ:def H matrix for L collective}.
\end{notationlist}

\subsection*{Grouped and hierarchical settings}

\begin{notationlist}
    \item[$D$, $S$, $N$]
    The coordinate set $D:=\{1,\cdots,d\}$, a partition $S=\{\alpha_1,\cdots,\alpha_N\}$ of $D$, and its number of groups; see the opening of \Cref{sec:multiple variables approach}.

    \item[$\alpha$, $\alpha^c$, $|\alpha|$]
    A group of coordinate indices, its complement in $D$, and its cardinality; see the opening of \Cref{sec:two variables approach}.

    \item[$\bfX_\alpha$, $\calX_\alpha$, $\mu_\alpha$]
    The subvector of $\bfX$ indexed by $\alpha$, its state space, and its probability distribution.  
    Independence of the groups is written $\mu=\bigotimes_{\alpha\in S}\mu_\alpha$; see \Cref{sec:two variables approach,sec:multiple variables approach}.

    \item[$\bfm=(m_\alpha)_{\alpha\in S}$]
    The tuple containing the number of features assigned to each group; see \Cref{sec:multiple variables approach}.

    \item[$g^\alpha$, $g=(g^\alpha)_{\alpha\in S}$]
    The feature map acting only on group $\alpha$ and the associated separated feature map on all groups; see \Cref{sec:multiple variables approach}.

    \item[$\calG_{m_\alpha}^\alpha$, $\calG_{\bfm}$]
    The admissible feature maps on group $\alpha$ and their product class $\calG_{\bfm}=\prod_{\alpha\in S}\calG_{m_\alpha}^\alpha$; see \Cref{sec:multiple variables approach}.

    \item[$\calF_{\bfm}^{bi}$, $\calF_{\bfm}^{mul}$]
    The classes of bilinear and multilinear profile functions, identified with matrices and order-$N$ tensors, respectively; see \Cref{subsec:bilinear regression function,subsec:multilinear regression function}.

    \item[$\calE_\alpha^{bi}(g)$, $\calE_S^{mul}(g)$]
    The best approximation errors for bilinear and multilinear profile functions; see \eqref{equ:approximation error bilinear} and \eqref{equ:approximation error multilinear}.

    \item[$\calE(g)$]
    The best grouped approximation error when the profile function is only required to be measurable; see \eqref{equ:def E}.

    \item[$\Sigma(g)$, $\Sigma_\alpha(g^\alpha)$]
    The closed subspaces of square-integrable functions measurable with respect to $g(\bfX)$ and $g^\alpha(\bfX_\alpha)$, respectively; see \eqref{equ:def Sigma g} and \eqref{equ:def Sigma beta g}.

    \item[$id^\alpha$]
    The identity feature map on the coordinates indexed by $\alpha$; see the paragraph preceding Proposition \ref{prop:approximation error decomposition multi}.

    \item[$\mathrm{rankf}_\alpha(v)$]
    The $\alpha$-feature-rank of $v$, namely the smallest number of features of $\bfx_\alpha$ needed to express $v$ as a function of those features and $\bfx_{\alpha^c}$; see the definition at the beginning of \Cref{sec:toward hierarchical formats}.
\end{notationlist}

\end{document}